\documentclass[11pt]{article}
\usepackage{amsmath}
\def\disp{\displaystyle}

\def\Limsup{\mathop{{\rm Lim}\,{\rm sup}}}

\def\tto{\;{\lower 1pt \hbox{$\rightarrow$}}\kern -10pt
\hbox{\raise 2pt \hbox{$\rightarrow$}}\;}
\def\Hat{\widehat}
\def\Tilde{\widetilde}
\def\tilde{\widetilde}
\def\Bar{\overline}
\def\ra{\rangle}
\def\la{\langle}
\def\ve{\varepsilon}
\def\B{I\!\!B}
\def\h{\hfill\Box}
\def\R{I\!\!R}
\def\N{I\!\!N}
\def\ox{\bar{x}}
\def\ow{\bar{w}}
\def\oy{\bar{y}}
\def\oz{\bar{z}}

\def\ot{\bar{t}}
\def\oq{\bar{q}}
\def\o{\omega}

\def\epi{\mbox{\rm epi}\,}
\def\dim{\mbox{\rm dim}\,}
\def\dom{\mbox{\rm dom}\,}

\def\i{\mbox{\rm int}\,}
\def\substack#1#2{{\scriptstyle{#1}\atop\scriptstyle{#2}}}
\def\subsubstack#1#2#3{{{{\scriptstyle{#1}}\atop{\scriptstyle{#2}}}
\atop{\scriptstyle{#3}}}}
\def\h{\hfill\triangle}
\def\dn{\downarrow}
\def\O{\Omega}
\def\ph{\varphi}
\def\emp{\emptyset}
\def\st{\stackrel}
\def\oR{\Bar{\R}}
\def\lm{\lambda}
\def\gg{\gamma}
\def\dd{\delta}
\def\al{\alpha}

\def\N{I\!\!N}

\newcounter{lk}

\begin{document}
\begin{center}
\vspace*{0.3in} {\bf SUBGRADIENTS OF MINIMAL TIME FUNCTIONS\\
UNDER MINIMAL REQUIREMENTS}\\[3ex]
BORIS S. MORDUKHOVICH\footnote{Department of Mathematics, Wayne
State University, Detroit, MI 48202, USA (email:
boris@math.wayne.edu). Research of this author was partially
supported by the US National Science Foundation under grant
DMS-0603846 and by the Australian Research Council under grant
DP-12092508.} and NGUYEN MAU NAM\footnote{Department of
Mathematics, The University of Texas--Pan American, Edinburg,
TX 78539--2999, USA (email: nguyenmn@utpa.edu).}\\[2ex]
\end{center}
\small This paper concerns the study of a broad class of minimal
time functions corresponding to control problems with constant convex
dynamics and closed target sets in arbitrary Banach spaces. In
contrast to other publications, we do not impose any nonempty
interior and/or calmness assumptions on the initial data and
deal with generally non-Lipschitzian minimal time functions. The
major results present refined formulas for computing various subgradients
of minimal time functions under minimal requirements in both cases of convex
and nonconvex targets. Our technique is based on advanced tools of variational
analysis and generalized differentiation.\\[1ex]
{\em Keywords:} Variational analysis and optimization, minimal
time functions, Minkowski gauges, generalized differentiation,
subdifferentials and normal cones\\[1ex]
{\em 2000 Mathematical Subject Classification:} 49J52, 49J53,
90C31
\newtheorem{Theorem}{Theorem}[section]
\newtheorem{Proposition}[Theorem]{Proposition}
\newtheorem{Remark}[Theorem]{Remark}
\newtheorem{Lemma}[Theorem]{Lemma}
\newtheorem{Corollary}[Theorem]{Corollary}
\newtheorem{Definition}[Theorem]{Definition}
\newtheorem{Example}[Theorem]{Example}
\renewcommand{\theequation}{\thesection.\arabic{equation}}
\normalsize

\section{Introduction}
\setcounter{equation}{0}

Consider the {\em minimal time problem} with constant dynamics
given by
\begin{eqnarray}\label{mtp}
\mbox{minimize }\;t\ge 0\;\mbox{ subject to }\;(x+t
F)\cap\O\ne\emp,\quad x\in X,
\end{eqnarray}
where $X$ is an arbitrary Banach space of state variables,
$\O\subset X$ is a closed {\em target} set, and $F\subset X$ is a
closed, bounded, and convex set describing the {\em constant
dynamics} $\dot x\in F$ to attain the target set $\O$ from the
state $x\in X$. We refer the reader to
\cite{bar,can,cowo1,heng,mn10,sova,wozh} and the bibliographies
therein for various results and discussions on the minimal time
problems and their applications, particularly to control and
optimization.

The main attention of this paper is paid to the optimal value
function in problem \eqref{mtp} known as the {\em minimal time
function} and defined by
\begin{equation}\label{mt}
T^F_\O(x):=\inf\big\{t\ge 0\big|\;\O\cap(x+tF)\ne\emp\big\}.
\end{equation}
The requirements on the initial data $(X,\O,F)$ of \eqref{mtp}
imposed above are our {\em standing assumptions} in this paper.
Observe that we do {\em not} assume the standard {\em interiority
condition} $0\in\i F$, which is a conventional requirement on $F$
in the study of the minimal time function \eqref{mt} ensuring, in
particular, the Lipschitz continuity of \eqref{mt} as well as of
the corresponding {\em Minkowski gauge}
\begin{eqnarray}\label{mk}
\rho_F(u):=\inf\big\{t\ge 0\big|\; u\in tF\big\},\quad u\in X,
\end{eqnarray}
generating \eqref{mt} under the interiority condition by
\begin{eqnarray}\label{mtg}
T^F_\O(x)=\disp\inf_{w\in\Omega}\rho_F( w-x),\quad x\in X,
\end{eqnarray}
where $\rho_F(u)=\inf\{t>0|\;t^{-1}u\in F\}$ in this case.
Representation \eqref{mtg} with the Lipschitz continuous gauge
\eqref{mk} relates the minimal time function $T^F_\O(x)$ to the
classical {\em distance function} of the set $\O$ defined by
\begin{equation}\label{df}
\mbox{dist}(x;\O):=\inf_{y\in\O}\|y-x\|,\quad x\in\Omega,
\end{equation}
which corresponds to \eqref{mt} when $F=\B$, the closed unit ball
in $X$. In fact, the vast majority of methods and results
developed in the study of the minimal time function \eqref{mtg}
under the interiority requirement $0\in\i F$ are inspired by their
counterparts for the distance function \eqref{df}; see more
details and discussions in the reference above. In the absence of
the latter requirement the minimal time function may be quite
different from the distance one; e.g., for $F=[0,1]\subset\R$ and
$\O=(-\infty,0]$ we have the expression
$$T^F_\O(x)=\begin{cases}
0&\mbox{if }\;x\in\O,\\
\infty&\text{otherwise.}
\end{cases}$$

It is worth noting that functions of type \eqref{mt} arise not only
in the control framework and have not only the ``minimal time"
interpretation. Their importance has been well recognized in
approximation theory; see, e.g., \cite{demy,ln}. Furthermore,
functions of type \eqref{mt} belong to a broader class of the
so-called {\em marginal functions}
\begin{eqnarray}\label{mar}
\mu(x):=\disp\inf_{ w\in\O(x)}\ph(x,\o),\quad x\in\O,
\end{eqnarray}
describing, in particular, {\em optimal values} in general
problems of parametric optimization and playing a significant role
in sensitivity, stability, and other aspects of variational
analysis and its applications; see, e.g.,
\cite{bz,mor06a,mor06b,mny,rw,t} and the references therein.
However, the {\em special structure} of the cost
function/Minkowski gauge in \eqref{mtg} is {\em crucial} for the
most interesting results obtained for the minimal time and
distance functions and cannot be deduced from those known for more
general classes of marginal functions \eqref{mar}.\vspace*{0.05in}

A characteristic feature of the minimal time function \eqref{mt}
is its {\em intrinsic nonsmoothness}, which requires the usage of
appropriate tools of generalized differentiation. A number of
results for evaluating various subdifferentials of \eqref{mt} were
given in \cite{cowo,cowo1,heng,mn10,wozh} under the underlying
assumption $0\in\i F$, which ensures that the Lipschitz continuous
function $T^F_\O(x)$ behaves similarly to the distance function
\eqref{df} from the viewpoint of generalized differentiation. It
is definitely not the case when the assumption $0\in\i F$ is
violated.

To the best of our knowledge, the first effort in dealing with the
minimal time functions of type \eqref{mt} in the absence of the
interiority condition $0\in\i F$ was made in \cite{jh}, where
certain formulas for evaluating their proximal and Fr\'echet
subdifferentials were obtained. However, the major results in the
out-of-set case $\ox\notin\O$ were derived in \cite{jh} under the
{\em calmness} property \cite{rw} of $T^F_\O(\cdot)$ at $\ox$
meaning that
\begin{eqnarray}\label{calm}
\big|T^F_\O(x)-T^F_\O(\ox)\big|\le\kappa\|x-\ox\|\;\;\mbox{ for
all $x$ near $\ox$}
\end{eqnarray}
with some constant $\kappa>0$, which is a ``one-point" refinement
of the classical Lipschitz continuity of the minimal time function
discussed above.\vspace*{0.05in}

The primary goal of this paper is to develop subdifferential
properties of the minimal time function \eqref{mt} with {\em no
imposing either} the interiority condition $0\in\i F$ {\em or} the
calmness condition \eqref{calm}. Besides the pure theoretical
interest of clarifying what is possible to get without the
aforementioned requirements, the major motivation for our study
comes from the application to the {\em generalized
Fermat-Torricelli problem} of finding a point at which the sum of
its distances to the given closed (convex and non convex) sets is
minimal. The latter problem is comprehensively studied in the
associated paper \cite{mn10a} from both qualitative and
quantitative/numerical viewpoints.

We pay the main attention to the two robust limiting constructions
by Mordukhovich: the {\em basic/limiting} and {\em singular}
subdifferentials for minimal time functions. The first of them was
studied in our recent paper \cite{mn10} in the case of $0\in\i F$
while the second one, being trivial for Lipschitzian functions,
was not considered in \cite{mn10} or anybody else in the literature
on minimal time functions. As a preliminary technical step (but of
its own interest) we evaluate $\ve$-subdifferentials of the Fr\'echet
type for \eqref{mt}. The latter construction reduces to the usual
Fr\'echet subdifferential studied in \cite{jh}, while we need its
$\ve$-enlargements in the general Banach space setting. Note that
some results obtained here for Fr\'echet subgradients of \eqref{mt}
recover those from \cite{jh}, while the most of them are new in the
settings under consideration, even in the case of convex data with
no calmness assumption.\vspace*{0.05in}

The rest of the paper is organized as follows. Section~2 contains
preliminaries from generalized differentiation used in the
formulations and proofs of the main results.

Section~3 concerns general (non-subdifferential) properties of
minimal time functions important for their own sake and useful for
the subsequent study of subdifferentials.

Section~4 deals with $\ve$-subdifferentials of \eqref{mt} at
$\ox\in X$ considering both in-set $\ox\in\O$ (easier) and
out-of-set $\ox\notin\O$ (more difficult) cases. The crucial
result in the latter case is representing $\ve$-subgradients of
the minimal time function via appropriate normals at perturbed
projections on the target with proofs based on variational
arguments.

In Sections~5--7 we present the main results of the paper related to
evaluating basic and singular subgradients of minimal time functions
in both convex and nonconvex settings. Most of the results obtained in
these lines are new even for the case of $0\in\i F$ and are illustrated
by numerical examples.

Section~5 is particularly devoted to upper estimates and precise
representations of the basic and singular subdifferentials of \eqref{mt}
at in-set points $\ox\in\O$ of general nonconvex target sets. It contains
upper estimates and equalities for evaluating basic and singular
subgradients of the minimal time function $T^F_\O$ via the limiting normals
to the target $\O$ and appropriate characteristics of the dynamics $F$.

Section~6 concerns upper estimates and equalities for the basic and singular
subdifferentials of $T^F_\O$ and their one-sided counterparts at out-of-set
points $\ox\notin\O$ of the general target set $\O$. We derive two types of
results in this direction: those expressed via limiting normals to $\O$ at
projection points and those involving the limiting normal cone to the
corresponding enlargements of the target.

Section~7 is devoted to the minimal time problem \eqref{mtp} with convex data.
The exact calculations of the convex subdifferential of \eqref{mt}
obtained here recover some results of \cite{jh} but without the calmness
condition and also provide new subdifferential formulas involving the
Minkowski gauge \eqref{mk} in the absence of the interiority condition
$0\in\i F$. Besides computing the convex subdifferential of \eqref{mt},
we give the exact evaluation of the singular subdifferential of the
convex minimal time function, which has never been consider in the
minimal time literature. It is worth also mentioning that the singular
subdifferential has not been systematically studied and applied in the
general framework of convex analysis.
\vspace*{0.05in}

Out notation is basically standard in variational analysis and
generalized differentiation; see, e.g., \cite{mor06a,rw}. Unless
otherwise stated, the space $X$ in question is arbitrary {\em
Banach} with the norm $\|\cdot\|$, the closed unit ball $\B$, and
the canonical pairing $\la\cdot,\cdot\ra$ between $X$ and its
topological dual $X^*$. As usual, the symbol $x_k\to\ox$ stands
for the norm convergence in $X$ while $x^*_k\st{w^*}{\to}x^*$ as
$k\in\N:=\{1,2,\ldots\}$ signifies the sequential weak$^*$
convergence in the dual space $X^*$. Given a set-valued mapping
$G\colon X\tto X^*$, we denote
\begin{eqnarray}\label{pk}
\begin{array}{ll}
\disp\Limsup_{x\to\ox}G(x):=\Big\{x^*\in X^*\Big|&\exists\;\mbox{
sequences
}\;x_k\to\ox,\;\;x^*_k\st{w^*}{\to}x^*\;\mbox{ as }\;k\to\infty\\
&\mbox{with }\;x^*_k\in G(x_k)\;\mbox{ for all }\;k\in\N\Big\}
\end{array}
\end{eqnarray}
the {\em sequential} Painlev\'e-Kuratowski upper/outer limit of
$G$ as $x\to\ox$. If no confusion arises, the symbol
$x\st{\O}{\to}x$ means that $x\to\ox$ with $x\in\O$ for a set
$\O$, while $x\st{\ph}{\to}\ox$ indicates that $x\to\ox$ with
$\ph(x)\to\ph(\ox)$ for an extended-real-valued function
$\ph\colon X\to\oR:=(-\infty,\infty]$.

\section{Preliminaries from Generalized Differentiation}
\setcounter{equation}{0}

Here we define the constructions of generalized differentiation in
variational analysis used in this paper and review some of their
properties. We mostly follow the book \cite{mor06a}, where the
reader can find comprehensive material in this direction with the
vast commentaries and references on these and related topics; cf.\
also \cite{bz,mor06b,rw,s} for additional issues.

Given a set $\O\subset X$ with $\ox\in\O$ and given $\ve\ge 0$,
the collection of {\em $\ve$-normals} to $\O$ at $\ox$ is
\begin{eqnarray}\label{1.2.1}
\Hat N_\ve(\ox;\O):=\disp\Big\{x^*\in
X^*\Big|\;\limsup_{x\st{\O}{\to}\ox}\frac{\la
x^*,x-\ox\ra}{\|x-\ox\|}\le\ve\Big\},\quad\ox\in\O,
\end{eqnarray}
with $\Hat N_\ve(\ox;\O)=\emp$ if $\ox\notin\O$ for convenience.
When $\ve=0$ in \eqref{1.2.1}, the set $\Hat N(\ox;\O):=\Hat
N_0(\ox;\O)$ is a cone known as the {\em Fr\'echet/regular normal
cone} to $\O$ at $\ox$. For {\em convex} sets $\O$ we have
\begin{eqnarray}\label{n-convex}
\Hat N_\ve(\ox;\O)=\big\{x^*\in X^*\big|\;\la
x^*,x-\ox\ra\le\ve\|x-\ox\|\;\mbox{ whenever
}\;x\in\O\big\},\quad\ox\in\O,
\end{eqnarray}
i.e., $\Hat N(\ox;\O)$ reduces to the normal cone of convex
analysis, while for nonconvex sets $\O$ the cone $\Hat N(\ox;\O)$
and its $\ve$-enlargement \eqref{1.2.1} do not generally possess
appropriate properties expected from natural notions of normals.
In particular, $\Hat N(\ox;\O)$ if often trivial $(=\{0\})$ for
boundary points of closed sets; there is no robustness and good
calculus for \eqref{1.2.1}, etc.

The situation dramatically changes when we consider the robust
sequential regularization \eqref{pk} of the set-valued mapping
\eqref{1.2.1} near $\ox$ defined by
\begin{eqnarray}\label{2.3.1}
N(\ox;\O):=\Limsup_\substack{x\to\ox}{\ve\dn 0}\Hat N_\ve(x;\O)
\end{eqnarray}
and known as the {\em basic/limiting/Mordukhovich normal cone} of
$\O$ at $\ox$. The latter cone enjoys a number of good properties
in the general Banach space setting and perfect ones in Asplund
spaces (including all reflexive) characterized as those Banach
spaces, where every separable subspace has a separable dual; see
\cite{bz,mor06a,p} for more details. In this paper we do not need
to impose the Asplund structure on $X$. Note that the normal cone
\eqref{2.3.1} and the corresponding subdifferentials are usually
{\em nonconvex} (in contrast to the majority of their known
counterparts), while their important properties and applications
are mainly based on {\em variational/extremal principles} of
variational analysis.

In this paper we employ the following three subgradient
constructions for extended-real-valued functions $\ph\colon
X\to\oR$ generated by normals \eqref{1.2.1} and \eqref{2.3.1} to
their epigraphs $\epi\ph:=\{(x,\mu)\in X\times\R|\;\mu\ge\ph(x)\}$.
For convenience we present these constructions in the equivalent
analytic forms. Given a function $\ph\colon X\to\oR$ and a point $\ox$
from its domain $\dom\ph:=\{x\in X|\;\ph(x)<\infty\}$,
the {\em $\ve$-subdifferential} of the Fr\'echet type of $\ph$ at
$\ox$ is given by
\begin{eqnarray}\label{2.3}
\Hat\partial_\ve\ph(\ox):=\Big\{x^*\in X^*\Big|\;\liminf_{
x\to\ox}\frac {\ph(x)-\ph(\ox)-\la x^*,x-\ox\ra}{\|x-\ox
\|}\ge-\ve\Big\},\quad\ve\ge 0,
\end{eqnarray}
with $\Hat\partial\ph(\ox):=\Hat\partial_0\ph(\ox)$.
For {\em convex} functions $\ph$ the $\ve$-subdifferential \eqref{2.3}
reduces to
\begin{eqnarray}\label{2.3a}
\Hat\partial_\ve\ph(\ox)=\big\{x^*\in X^*\big|\;\la x^*,x-\ox\ra\le
\ph(x)-\ph(\ox)+\ve\|x-\ox\|\;\mbox{ whenever }\;x\in X\big\}.
\end{eqnarray}

The {\em basic subdifferential} $\partial\ph(\ox)$ and {\em singular
subdifferential} $\partial^\infty\ph(\ox)$ of Mordukhovich are
generated, respectively, by ``slant" and ``horizontal" normals to
$\epi\ph$ at $(\ox,\ph(\ox))$ in the sense of \eqref{2.3.1} and
can be defined analytically as
\begin{eqnarray}\label{2.4}
\partial\ph(\bar x):=\Limsup_\substack{x\xrightarrow{\ph}\bar x}{\ve\dn 0}
\Hat\partial_\ve\ph(x),
\end{eqnarray}
\begin{eqnarray}\label{2.4.1}
\partial^\infty\ph(\bar x):=\Limsup_\subsubstack{x\xrightarrow{\ph}\bar x}
{\lm\dn 0}{\ve\dn 0}\lm\Hat\partial_\ve\ph(x).
\end{eqnarray}
It is worth observing (although it is not used in the paper) that
we can equivalently put $\ve=0$ in \eqref{2.4} and \eqref{2.4.1}
if $X$ is Asplund and $\ph$ is lower semicontinuous (l.s.c.)
around $\ox$.

Recall that the Fr\'echet subdifferential $\Hat\partial\ph(\ox)$
reduces to the classical Fr\'echet derivative of $\ph$ at $\ox$ if
$\ph$ is Fr\'echet differentiable at $\ox$, while the basic
subdifferential \eqref{2.4} reduces to the classical derivative
$\partial\ph(\ox)=\{\nabla\ph(\ox)\}$ if $\ph$ is strictly
differentiable at $\ox$ in the sense of
\begin{eqnarray*}
\disp\lim_\substack{x\to\ox}{u\to\ox}\frac{\ph(x)-\ph(u)-\la\nabla\ph(\ox),
x-u\ra}{\|x-u\|}=0,
\end{eqnarray*}
which is automatic when $\ph$ is $C^1$ around $\ox$. If $\ph$ is
convex, both $\Hat\partial\ph(\ox)$ and $\partial\ph(\ox)$ agree
with the subdifferential of convex analysis.

For the singular subdifferential \eqref{2.4.1} we have
$\partial^\infty\ph(\ox)=\{0\}$ if $\ph$ is locally Lipschitzian
around $\ox$ in arbitrary Banach spaces. In fact, the latter singular
subdifferential condition is a full {\em characterization} of the
local Lipschitzian property under some additional assumptions,
which are automatic in finite dimensions; see \cite[Theorem~3.52]{mor06a}.
Thus the singular subdifferential carries nontrivial information only
for non-Lipschitzian functions, which is not the case for the minimal
time function \eqref{mt} under the interiority condition $0\in\i F$.

\section{General Properties of Minimal Time Functions}
\setcounter{equation}{0}

In this section we collect some properties of the minimal time
function \eqref{mt}, which are not related to generalized
differentiation. They are of their own interest while most of them
are widely used in the subsequent sections for deriving
subdifferential results of the paper. Note that, under our {\em
standing assumptions} made in Section~1 and imposed in what
follows, the minimal time function is merely extended-real-valued
$T^F_\O\colon X\to\oR$ and does not share many common properties with
the distance function \eqref{df} as in the case of $0\in\i
F$.\vspace*{0.05in}

For the given target set $\O$, consider the family of its {\em
enlargements}
\begin{equation}\label{enl}
\O_r:=\big\{x\in X\big|\;T^F_\O(x)\le r\big\},\quad r>0,
\end{equation}
and establish the following relationship between $T^F_{\O_r}$ and
$T^F_\O$.

\begin{Proposition}\label{distance estimate} {\bf (minimal time functions
for enlargements of target sets).} Let $x\notin\O_r$ be such that
$T^F_\O(x)<\infty$. Then
\begin{equation}\label{enl1}
T^F_\O(x)=r+T^F_{\O_r}(x)\;\mbox{ whenever }\;r>0.
\end{equation}
\end{Proposition}
{\bf Proof.} Since $\O\subset\O_r$, we have
$t_1:=T^F_{\O_r}(x)<\infty$. By the definition of $T^F_{\O_r}(x)$,
for any $\ve>0$ there are $w_1\in\O_r$ and
$t_1\le\gamma_1<t_1+\ve$ satisfying
\begin{equation*}
w_1\in\O_r\cap(x+\gamma_1 F).
\end{equation*}
Then $T^F_\O(w_1)\le r$, and hence there are $w_2\in\O$ and
$\gamma_2<r+\ve$ such that
\begin{equation*}
 w_2\in\Omega\cap (w_1+\gamma_2F).
\end{equation*}
Consequently we get $w_2\in\O\cap(x+(\gamma_1+\gamma_2)F)$ by the
convexity of $F$. This gives
\begin{equation*}
T^F_\O(x)\le\gamma_1+\gamma_2\le T^F_{\O_r}(x)+r+2\ve,
\end{equation*}
which imply in turn that $T^F_\O(x)\le T^F_{\O_r}(x)+r$ due to the
arbitrary choice of $\ve>0$.

To justify the opposite inequality in \eqref{enl1}, denote
$t:=T^F_\O(x)>r$. Then for any $\ve>0$ there exist $\gamma$ with
$t\le\gamma<t+\ve$ and $w\in X$ satisfying the relationship
\begin{equation*}
 w\in\O\cap (x+\gamma F).
\end{equation*}
The above element $w\in\O$ can be represented as $w=x+\gamma q$
with some $q\in F$. Define further $w_r:=x+(\gamma-r)q$ and get
$w_r\in\O_r$ by $w\in\O\cap( w_r+rF)\ne\emptyset$. Thus
$w_r\in\O_r\cap(x+(\gamma-r)F)$, which implies the inequalities
\begin{equation*}
T^F_{\O_r}(x)\le\gamma-r\le T^F_\O(x)+\ve-r.
\end{equation*}
We therefore arrive at $T^F_{\O_r}(x)+r\le T^F_\O(x)$ and complete
the proof of the proposition. $\h$\vspace*{0.05in}

The next property is elementary while useful in what follows.

\begin{Proposition}\label{prop2} {\bf (minimal time functions with shifted arguments).}
For any $x\in\O_r$ with $r>0$ and any $t\ge 0$ we have
\begin{equation*}
T^F_\O(x-tq)\le r+t\;\mbox{ whenever }\;q\in F.
\end{equation*}
\end{Proposition}
{\bf Proof.} Fix $(x,r,t,q)$ in the formulation of the theorem and
and denote $\lm:=T^F_\O(x)$. Picking any $\ve>0$ and observing
that $\lm\le r$, find $\gg>0$ such that $\lm\le\gg<\lm+\ve$ and
$w\in X$ satisfying
\begin{equation*}
w\in\O\cap(x+\gamma F).
\end{equation*}
The latter directly implies the inclusions
\begin{equation*}
w\in\O\cap(x-tq+tq+\gg F)\subset\O\cap(x-tq+tF+\gg
F)\subset\O\cap\big(x-tq+(t+\gg)F\big).
\end{equation*}
It follows then that $T^F_\O(x-tq)\le\gamma+t\le t+\lm+\ve\le
t+r+\ve$, and hence $T^F_\O(x-tq)\le r+t$ by the arbitrary choice
of $\ve>0$. $\h$\vspace*{0.05in}

Now we justify an important result ensuring the representation
\eqref{mtg} of the minimal time function \eqref{mt} via the
Minkowski gauge \eqref{mk} with no interiority requirement $0\in\i
F$.

\begin{Proposition}\label{mtg1} {\bf (relationship between minimal
time and Minkowski functions).} Under the standing assumptions
made we have the representation
\begin{eqnarray*}
T^F_\O(x)=\inf_{w\in\O}\rho_F( w-x)\;\mbox{ for all }\;x\in X.
\end{eqnarray*}
\end{Proposition}
{\bf Proof.} Let us first show that $T^F_\O(x)=\infty$ if and only
if
\begin{equation}\label{infty}
\inf_{w\in\O}\rho_F( w-x)=\infty,\quad x\in X.
\end{equation}
Indeed, it follows from definition \eqref{mt} that
$T^F_\O(x)=\infty$ for some fixed $x\in X$ if and only if
$\O\cap(x+tF)=\emp$ whenever $t\ge 0$. The latter is equivalent to
the fact that
\begin{equation*}
\big\{t\ge 0\big|\;w-x\in tF\big\}=\emp\;\mbox{ for any }\;w\in\O,
\end{equation*}
which is the same as $\rho_F(w-x)=\infty$ for all $w\in\O$, i.e., \eqref{infty} holds.

Suppose now that $T^F_\O(x)<\infty$ and thus
$\inf_{w\in\O}\rho_F(w-x)<\infty$ for some fixed $x\in X$. Then for any $t\ge 0$ with
$\O\cap(x+tF)\ne\emp$ there is $w\in\O$ satisfying $w-x\in tF$, and hence $\rho_F(w-x)\le t$.
The latter implies that
\begin{equation*}
\inf_{w\in\O}\rho_F(w-x)\le t,
\end{equation*}
and so $\inf_{w\in\O}\rho_F(w-x)\le T^F_\O(x)$. Put further
$\gg:=\inf_{w\in\O}\rho_F(w-x)<\infty$ and, given $\ve>0$, find
$w\in\O$ satisfying
\begin{equation*}
\rho_F(w-x)<\gg+\ve.
\end{equation*}
Then there is $t\ge 0$ such that $t<\gg+\ve$ and $w-x\in tF$. This
implies that
\begin{equation*}
T^F_\O(x)\le t\le\gg+\ve,
\end{equation*}
and hence $T^F_\O(x)\le\gg=\inf_{w\in\O}\rho_F(w-x)$, which
completes the proof. $\h$\vspace*{0.05in}

Given $\ox\in X$ with $T^F_\O(\ox)<\infty$, consider the
(generalized, minimal time) {\em projection} of $\ox$ on the
target set $\O$ defined by
\begin{equation}\label{projection}
\Pi_\O^F(\ox):=\big(\ox+T^F_\O(\ox)F\big)\cap\O.
\end{equation}
It is not hard to check that if $\O$ is a compact set, the
projection $\Pi_\O^F(\ox)$ is always nonempty with $T^F_\O(\ox)=0$
if and only if $\ox\in\O$.\vspace*{0.05in}

The next result reveals a kind of linearity of the minimal time
functions with respect to projection points on arbitrary target
sets.

\begin{Proposition}\label{lin} {\bf (minimal time linearity with
respect to projections).}  Let $\ox\notin\O$, and let $\bar
w\in\Pi^F_\O(\ox)$. Then for any $\lambda\in(0,1)$ we have
\begin{equation}\label{lin1}
T^F_\O\big(\lambda\bar w+(1-\lambda)\bar
x\big)=(1-\lambda)T^F_\O(\ox).
\end{equation}
\end{Proposition}
{\bf Proof.} It follows that $\bar w\in\ox+tF$ for
$t:=T^F_\O(\ox)<\infty$. Then
\begin{equation*}
\lambda\bar w+(1-\lambda)\bar x=\bar  w+(1-\lambda)(\bar x-\bar
 w)\in\bar w-(1-\lambda)tF,
\end{equation*}
which implies the inclusion
\begin{equation*}
\bar w\in\O\cap\big(\lambda\bar w+(1-\lambda)\bar
x+(1-\lambda)tF\big),\quad 0<\lm<1.
\end{equation*}
Hence $T^F_\O(\lambda\bar w+(1-\lambda)\bar
x)\le(1-\lambda)t=(1-\lambda)T^F_\O(\ox)$ for such $\lm$, which
justifies the inequality ``$\le$" in \eqref{lin1}. To prove the
opposite inequality, denote $t_\lambda:=T^F_\O(\lambda\bar
w+(1-\lambda)\ox)<\infty$ and for any $\ve>0$ find $t_\lambda\le
\gamma<t_\lambda+\ve$ with
\begin{equation*}
\O\cap\big(\ox+\lambda(\bar w-\ox)+\gamma F\big)\ne\emp.
\end{equation*}
Thus we have that
\begin{equation*}
\O\cap\big(\ox+(\lambda t+\gamma)F\big)\ne\emp,
\end{equation*}
and so $T^F_\O(\ox)\le\lambda t+\gamma\le\lambda
T^F_\O(\ox)+t_\lambda+\ve$. It follows finally that
\begin{equation*}
(1-\lambda)T^F_\O(\ox)\le t_\lambda+\ve,
\end{equation*}
which completes the proof by passing to the limit as $\ve\dn 0$.
$\h$\vspace*{0.05in}

Let us now show that, not being Lipschitzian or calm under our
assumptions, the minimal time function \eqref{mt} enjoys the
desired {\em lower semicontinuity property} provided some
additional requirements needed for our subsequent applications.
Recall that the lower semicontinuity of an extended-real-valued
function $\ph\colon X\to\oR$ is equivalent to the closedness of
its {\em level sets} $\{x\in X|\;\ph(x)\le\al\}$ for all
$\al\in\R$.

\begin{Proposition}\label{lsc} {\bf (lower semicontinuity of minimal
time functions).} In addition to our standing assumptions, suppose
that the space $X$ is either finite-dimensional, or it is
reflexive and the target set $\O$ is convex. Then the minimal time
function \eqref{mt} is lower semicontinuous on its domain.
\end{Proposition}
{\bf Proof.} Fix any $\al\ge 0$ and show that the level set
\begin{equation*}
\mathcal{L}_\al:=\big\{x\in X\big|\;T^F_\O(x)\le\al\big\}
\end{equation*}
is closed under the assumptions made. Take an arbitrary sequence
$\{x_k\}\subset\mathcal{L}_\al$ with $x_k\to\ox$ as $k\to\infty$.
Then we have from $T^F_\O(x_k)\le\al$ and definition \eqref{mt}
that for every $k\in\N$ there is $t_k$ such that $0\le
t_k<\al+1/k$ and
\begin{equation*}
\O\cap(x_k+t_k F)\ne\emp,\quad k\in\N.
\end{equation*}
Fixing further $w_k\in\O$ with $w_k\in x_k+t_kF$, we find $q_k\in
F$ satisfying $w_k=x_k+t_k q_k$ for all $k\in\N$. Observe that the
sequences $\{t_k\}$ and $\{q_k\}$ are bounded. If $X$ is {\em
finite-dimensional}, we get without loss of generality that
$t_k\to\ot$ and $q_k\to\oq$ as $k\to\infty$ for some elements
$\ot\in[0,\al]$ and $\oq\in F$. Then $w_k=x_k+t_kq_k\to\ox+\bar
t\bar q\in\O$, and thus $T^F_\O(\ox)\le\bar t\le\al$.

If $X$ is {\em reflexive}, we may assume that $q_k$ converges
weakly to some $\oq$. It follows from the classical Mazur theorem
that a convex combination of elements from the sequence $\{q_k\}$
converge to $\oq$ strongly in $X$. By the closedness and convexity
of $F$ we conclude that $\oq\in F$, and the same properties of
$\O$ imply that $\ox+\ot\oq\in\O$. Thus $T^F_\O(\ox)\le\bar
t\le\al$ in this case as well, which completes the proof of the
proposition. $\h$\vspace*{0.05in}

Next we characterize the {\em convexity} property of the minimal
time function $T^F_\O(x)$.

\begin{Proposition}\label{conv} {\bf (convexity of minimal time functions).}
 The minimal time function \eqref{mt} is convex if and only if its target set
 $\O$ is convex.
\end{Proposition}
{\bf Proof.} Suppose that the target set $\O$ is convex and show
that in this case for any $x_1,x_2\in X$ and for any
$\lambda\in(0,1)$ we have
\begin{equation}\label{convexinequality}
T_\O^F\big(\lambda x_1+(1-\lambda)x_2\big)\le\lambda
T^\O_F(x_1)+(1-\lambda)T^F_\O(x_2).
\end{equation}
Since (\ref{convexinequality}) obviously holds if
$T^\O_F(x_1)=\infty$ or $T^\O_F(x_2)=\infty$, assume in what
follows that $t_1:=T^\O_F(x_1)<\infty$ and
$t_2:=T^\O_F(x_2)<\infty$. Then for any $\ve>0$ there are numbers
$\gamma_i$ with
\begin{equation*}
t_i\le\gamma_i<t_i+\ve\;\mbox{ and }\;\O\cap(x_i+\gamma_iF)\ne
\emp,\quad i=1,2.
\end{equation*}
Take $w_i\in\O\cap(x_i+\gamma_iF)$ and by the convexity of $\O$
and $F$ get $\lambda w_1+(1-\lambda) w_2\in\O$ and
\begin{align*}
\lambda w_1+(1-\lambda) w_2&\in \lambda
x_1+(1-\lambda)x_2+\lambda \gamma_1F+(1-\lambda)\gamma_2F\\
&\subset \lambda w_1+(1-\lambda)
w_2+\big(\lambda\gamma_1+(1-\lambda\big)\gamma_2)F.
\end{align*}
The latter implies the inequalities
\begin{align*}
T^F_\O\big(\lambda x_1+(1-\lambda)x_2\big)&\le\lambda
\gamma_1+(1-\lambda)\gamma_2\\
&\le\lambda T^F_\O(x_1)+(1-\lambda)T^F_\O(x_2)+\ve,
\end{align*}
which in turn justify (\ref{convexinequality}) by the arbitrary
choice of $\ve>0$.

To prove the converse statement of the proposition, observe that
\begin{equation*}
\O=\big\{x\in X\big|\;T^F_\O(x)\le 0\big\},
\end{equation*}
and thus $\O$ is obviously convex provided that $T^F_\O$ has this
property. $\h$\vspace*{0.05in}

The last result of this section establishes sufficient conditions
for {\em concavity} property of the minimal time function under
consideration.

\begin{Proposition}\label{conc} {\bf(concavity of minimal time functions).}
Assume that the complement $\O^c:=X\setminus\O$ of the target is
convex. Then the minimal time function \eqref{mt} is concave on
$\O^c$ provided that it is finite on this set.
\end{Proposition}
{\bf Proof.} If $T^F_\O$ is not concave on $\O^c$, then there are
$x_1,x_2\in\O^c$ and $0<\lambda<1$ such that
\begin{equation}\label{conc1}
T^F_\O\big(\lambda x_1+(1-\lambda)x_2\big)<\lambda
T^F_\O(x_1)+(1-\lambda)T^F_\O(x_2)<\infty.
\end{equation}
By definition \eqref{mt}, find $t<\lambda
T^F_\O(x_1)+(1-\lambda)T^F_\O(x_2)$ and $w\in\O$ satisfying
\begin{equation*}
w-(\lambda x_1+(1-\lambda)x_2)=tq
\end{equation*}
for some $q\in F$. Consider the points
\begin{equation*}
u_i:=x_i+\dfrac{tq}{\lambda
T^F_\O(x_1)+(1-\lambda)T^F_\O(x_2)}T^F_\O(x_i),\quad i=1,2,
\end{equation*}
and observe that $u_1,u_2\in\O^c$. Indeed, assuming for
definiteness that $u_1\in\O$ yields that
\begin{equation*}
T^F_\O(x_1)\le\dfrac{t T^F_\O(x_1)}{\lambda
T^F_\O(x_1)+(1-\lambda)T^F_\O(x_2)}<T^F_\O(x_1),
\end{equation*}
a contradiction. At the same time we have the inclusion $w=\lambda
u_1+(1-\lambda)u_2\in\O$, which is impossible due to the convexity
of $\O^c$. Combining all the above shows that condition
\eqref{conc1} does not hold under the assumptions made, and thus
$T^F_\O$ is concave on $\O^c$. $\h$

\section{$\ve$-Subgradients of Minimal Time Functions}
\setcounter{equation}{0}

This section is devoted to evaluating $\ve$-subgradients
\eqref{2.3} of the minimal time function \eqref{mt} as $\ve\ge 0$
via corresponding characteristics of the target and dynamics sets
therein at both in-set and out-of-set points of the target in the
general Banach space setting. In particular, our results for
$\ve=0$ provide evaluations of Fr\'echet subgradients of
\eqref{mt} with no interiority and/or calmness assumptions
essentially used in previous methods and results for this
case.\vspace*{0.05in}

We first consider {\em in-set} points $\ox\in\O$. Involving the
{\em support function} of the dynamics
\begin{equation}\label{sf}
\sigma_F(x^*):=\sup_{x\in F}\la x^*,x\ra,\quad x^*\in X^*,
\end{equation}
and the exact {\em dynamics bound}
\begin{eqnarray}\label{db}
\|F\|:=\sup\big\{\|q\|\;\mbox{ over }\;q\in F\big\},
\end{eqnarray}
define the following {\em support level set}:
\begin{equation}\label{C1}
C^*_\ve:=\big\{x^*\in X^*\big|\;\sigma_F(-x^*)\le
1+\ve\|F\|\big\},\quad\ve\ge 0,
\end{equation}
which is denoted by $C^*$ if $\ve=0$. Let us begin with upper
estimating the $\ve$-subdifferential of \eqref{mt} via the support
set \eqref{C1} of the dynamics and the set of $\ve$-normals
\eqref{1.2.1} to the target.

\begin{Proposition}\label{Frechet1} {\bf (upper estimate of
$\ve$-subdifferentials of minimal time functions at in-set
points).} Let $\ox\in\O$. Then we have
\begin{equation*}
\Hat\partial_\ve T^F_\O(\ox)\subset\Hat N_\ve(\ox;\O)\cap
C^*_\ve\;\mbox{ for any }\;\ve\ge 0.
\end{equation*}
\end{Proposition}
{\bf Proof.} Fix an arbitrary subgradient $x^*\in\Hat\partial_\ve
T^F_\O(\ox)$. By definition \eqref{2.3} of the
$\ve$-subdifferential, for any $\eta>0$ find $\delta>0$ such that
\begin{align*}
\la x^*, x-\ox\ra &\le T^F_\O(x)-T^F_\O(\ox)+(\ve+\eta)\|x-\ox\|\\
&\le T^F_\O(x)+(\ve+\eta)\|x-\ox|\|
\end{align*}
whenever $x\in\ox+\dd\B$; this takes into account that
$T^F_\O(\ox)=0$ on $\O$. It follows that
\begin{align*}
\la x^*,x-\ox\ra\le(\ve+\eta)\|x-\ox\|\;\mbox{ for all }\;x\in\O,
\end{align*}
and thus $x^*\in\Hat N_\ve(\ox;\O)$. Fix further any $q\in F$ and
get
\begin{align*}
\la x^*,-tq\ra &\le T^F_\O(\ox-tq)+(\ve+\eta)\|tq\|\\
&\le t+t(\ve+\eta)\|F\|
\end{align*}
when $t>0$ is sufficiently small. Since $\eta>0$ is also
arbitrarily small, the latter implies that $\sigma_F(-x^*)\le
1+\ve\|F\|$ and completes the proof of the proposition.
$\h$\vspace*{0.05in}

The next result provides a certain opposite estimate to
Proposition~\ref{Frechet1}.

\begin{Proposition}\label{Frechet2} {\bf (lower estimate of
$\ve$-subdifferentials of minimal time functions at in-set
points).} Let $\ox\in\O$, and let $\ve\ge 0$. Then for any $x^*\in
\Hat N_\ve(\ox;\O)\cap C^*_\ve$ we have
\begin{equation}\label{fr}
x^*\in\Hat\partial_{\mu\ve}T^F_\O(\ox)\;\mbox{ with
}\;\mu=\mu(x^*):=1+2\|F\|\cdot\|x^*\|.
\end{equation}
\end{Proposition}
{\bf Proof.} Arguing by contradiction, suppose that $x^*\notin
\Hat\partial_{\mu\ve}T^F_\O(\ox)$. Then
\begin{equation*}
\liminf_{x\to\ox}\dfrac{T^F_\O(x)-T^F_\O(\ox)-\la x^*,
x-\ox\ra}{||x-\ox||}<-\mu\ve,
\end{equation*}
and thus there exist $\gamma>0$ and a sequence $x_k\to\ox$ such
that
\begin{equation*}
T^F_\O(x_k)-\la x^*,x_k-\ox\ra\le(-\mu\ve-\gamma)\|x_k-\ox\|,\quad
k\in\N.
\end{equation*}
It follows that $x_k\notin\O$ for $k$ sufficiently large, since
otherwise it contradicts the fact that $x^*\in\Hat N_\ve(\ox;\O)$
due to $\mu\ve+\gamma>\ve$. This also implies for such $k$ that
\begin{equation*}
0<T^F_\O(x_k)\le\|x^*\|\cdot\|x_k-\ox\|,
\end{equation*}
and hence $T^F_\O(x_k)\to 0$ as $k\to\infty$. Since
$\|x_k-\ox\|^2>0$, for each $k$ sufficiently large there are
$t_k\ge 0$, $w_k\in\O$, and $q_k\in F$ satisfying
\begin{equation*}
w_k=x_k+t_kq_k\;\mbox{ and }\;T^F_\O(x_k)\le
t_k<T^F_\O(x_k)+\|x_k-\ox\|^2.
\end{equation*}
Consequently we have the relationships
\begin{align*}
\la x^*,w_k-\ox\ra&=\la x^*,x_k-\ox\ra+t_k\la x^*,
q_k\ra\\
&\ge\la x^*,x_k-\ox\ra+t_k(-1-\ve\|F\|)\\
&\ge\la x^*,x_k-\ox\ra+(T^F_\O(x_k)+\|x_k-\ox\|^2)(-1-\ve\|F\|)\\
&=\la x^*,x_k-\ox\ra-T^F_\O(x_k)-(1+\ve\|F\|)||x_k-\ox||^2-\ve T^F_\O(x_k)\|F\|\\
&\ge(\mu\ve+\gamma-\ve\|x^*\|\cdot\|F\|)\|x_k-\ox\|-(1+\ve\|F\|)\|x_k-\ox\|^2.
\end{align*}
On the other hands, it follows from $w_k\xrightarrow{\O}\ox$ and
$x^*\in\Hat N_\ve(\ox;\O)$ that
\begin{equation*}
\la x^*,w_k-\ox\ra\le(\ve+\nu)\|w_k-\ox\|
\end{equation*}
for any $\nu>0$ and $k$ sufficiently large. Observe also that
\begin{equation*}
\|w_k-\ox\|\le\|x_k-\ox\|+t_k\|F\|\le(1+\|x^*\|\cdot
\|F\|)\|x_k-\ox\|+\|x_k-\ox\|^2\|F\|.
\end{equation*}
Comparing these inequalities and letting $\nu\dn 0$ and
$k\to\infty$, we get the estimate
\begin{equation*}
\mu\ve+\gamma-\ve\|x^*\|\cdot\|F\|\le\ve(1+\|x^*\|\cdot\|F\|)
\end{equation*}
Taking into account the definition of $\mu$ in \eqref{fr}, we
arrive at a contradiction and thus complete the proof of the
proposition. $\h$\vspace*{0.05in}

Let us now turn to the {\em out-of-set} case of $\ox\notin\O$. The
following important result is an extension of
\cite[Theorem~3.5]{mn10} established under the interiority
assumption $0\in\i F$. The proof is based on variational arguments
involving the seminal {\em Ekeland variational principle}.

\begin{Theorem}\label{out-of-set Frechet} {\bf ($\ve$-subgradients
of minimal time functions at out-of-set points via perturbed
normals to target sets).} Let $\ox\notin\O$ with
$T^F_\O(\ox)<\infty$. Then for every $x^*\in\Hat\partial_\ve
T^F_\O(\ox)$, $\ve\ge 0$, and $\eta>0$ there is $\ow\in\O$
satisfying the relationships
\begin{eqnarray}\label{e0}
x^*\in\Hat N_{\ve+\eta}(\ow;\O)\;\mbox{ and }\;\|\ox-\ow\|\le\|F\|
T^F_\O(\ox)+\eta.
\end{eqnarray}
\end{Theorem}
{\bf Proof.} Fix $(x^*,\ve,\eta)$ from the formulation of the
theorem and, using the $\ve$-subdifferential construction
\eqref{2.3}, find $\dd>0$ such that
\begin{eqnarray}\label{e1}
\la x^*, x-\ox\ra\le
T^F_\O(x)-T^F_\O(\ox)+\Big(\ve+\dfrac{\eta}{2}\Big)\|x-\ox\|\;\mbox{
for all }\;x\in\ox+\dd\B.
\end{eqnarray}
The minimal time definition \eqref{mt} ensures the existence of
$t\ge 0$, $\Tilde w\in\O$, and $q\in F$ satisfying
\begin{eqnarray}\label{e2}
T^F_\O(\ox)\le t<T^F_\O(\ox)+\Tilde\eta^2\;\mbox{ and }\;\Tilde
w=\ox+tq\;\mbox{ with
}\;\Tilde{\eta}:=\min\Big\{\dfrac{\delta}{2},\dfrac{\eta}{2+\|F\|},1\Big\}.
\end{eqnarray}
It follows from \eqref{e1} and \eqref{e2} that for any
$w\in\O\cap(\Tilde w+\dd\B)$ we have the estimates
\begin{align*}
\la x^*,w-\Tilde w\ra&\le T^F_\O(w-\Tilde w+\ox)-T^F_\O(\ox;\O)+
\Big(\ve+\dfrac{\eta}{2}\Big)\|w-\Tilde w\|\\
&\le T^F_\O(w-tf)-T^F_\O(\ox)+\Big(\ve+\dfrac{\eta}{2}\Big)\|w-\Tilde w\|\\
&\le t-T^F_\O(\ox)+\Big(\ve+\dfrac{\eta}{2}\Big)\|w-\Tilde w\|\\
&\le\Big(\ve+\dfrac{\eta}{2}\Big)\|w-\Tilde w\|+\Tilde\eta^2.
\end{align*}
Consider further the complete metric space $E:=\O\cap(\Tilde
w+\dd\B)$ and define a continuous function $\ph\colon E\to\R$ on
it by
\begin{eqnarray}\label{ph}
\ph(w):=-\la x^*,w-\Tilde
w\ra+\Big(\ve+\dfrac{\eta}{2}\Big)\|w-\Tilde w\|+\Tilde\eta^2,\quad
w\in E.
\end{eqnarray}
We conclude from the constructions and estimates above that
\begin{eqnarray*}
\ph(\Tilde w)\le\inf_{w\in E}\ph(w)+\Tilde\eta^2.
\end{eqnarray*}
Applying the Ekeland variational principle to $\ph$ on $E$ allows
us to find $\ow\in E$ such that
\begin{eqnarray*}
\|\Tilde w-\ow\|<\Tilde\eta\;\mbox{ and
}\;\ph(\ow)\le\ph(w)+\Tilde\eta\|w-\ow\|\;\mbox{ whenever }\;w\in
E.
\end{eqnarray*}
This means by the definition of $\ph$ in \eqref{ph} that
\begin{align*}
-\la x^*,\ow-\Tilde w\ra+\Big(\ve+\dfrac{\eta}{2}\Big)\|\ow-\Tilde
w\| +\Tilde\eta^2\le-\la x^*,w-\Tilde
w\ra+\Big(\ve+\dfrac{\eta}{2}\Big)\|w-\Tilde
w\|+\Tilde\eta^2+\Tilde\eta\|w-\ow\|
\end{align*}
for all $w\in E$. Taking into account the construction of
$\Tilde\eta$ in \eqref{e2}, we get
\begin{eqnarray}\label{e3}
\la
x^*,w-\ow\ra\le\Big(\ve+\dfrac{\eta}{2}+\Tilde\eta\Big)\|w-\ow\|
\le(\ve+\eta)\|w-\ow\|.
\end{eqnarray}
It follows furthermore that
\begin{eqnarray*}
\|w-\Tilde w\|\le\|w-\ow\|+\|\ow-\Tilde
w\|<2\Tilde\eta<\delta\;\mbox{ for any
}\;w\in\O\cap(\ow+\Tilde\eta\B).
\end{eqnarray*}
This ensures that $\O\cap(\ow+\Tilde\eta\B)\subset E$ and hence
$x^*\in\Hat N_{\ve+\eta}(\ow;\O)$ by \eqref{1.2.1} and \eqref{e3}.
Employing finally the choice of $(t,q,\Tilde w,\Tilde\eta)$ in
\eqref{e2}, we get
\begin{align*}
\|\ox-\ow\|&\le\|\ox-\Tilde w\|+\|\Tilde w-\ow\|\le t\|q\|+\Tilde\eta\\
&\le\|F\|\big(T^F_\O(\ox)+\Tilde\eta^2\big)+\Tilde\eta\le\|F\|T^F_\O(\ox)+
\Tilde\eta(\|F\|+1)\\
&\le\|F\|T^F_\O(\ox)+\eta,
\end{align*}
which justifies the remaining estimate in \eqref{e0} and completes
the proof of theorem. $\h$\vspace*{0.05in}

Next result fully describes behavior of the support function
\eqref{sf} at $\ve$-subgradients of the minimal time function
\eqref{mt} taken at $\ox\notin\O$ via the dynamics bound
\eqref{db}.

\begin{Proposition}\label{norm estimate} {\bf (relationship between
dynamics and $\ve$-subgradients of minimal time functions at
out-of-set points).} Let $\ox\notin \O$ and $T^F_\O(\ox)<\infty$
for \eqref{mt}. Then for any $x^*\in\Hat\partial_{\ve}T^F_\O(\ox)$
we have the two-sided estimates
\begin{equation}\label{supp}
1-\ve\|F\|\le\sigma_F(-x^*)\le 1+\ve\|F\|,\quad\ve\ge 0.
\end{equation}
\end{Proposition}
{\bf Proof.} Fix $\ve\ge 0$ and
$x^*\in\Hat\partial_{\ve}T^F_\O(\ox)$. Picking an arbitrary number
$\gg>0$ and using the $\ve$-subgradient definition \eqref{2.3},
find $\delta>0$ such that
\begin{equation*}
\la x^*,x-\ox\ra\le
T^F_\O(x)-T^F_\O(\ox)+(\ve+\gamma)\|x-\ox\|\;\mbox{ for all
}\;x\in\ox+\ve\B.
\end{equation*}
Let $r:=T^F_\O(\ox)$, which ensures that $\ox$ belongs to the
enlargement $\O_r$ defined in \eqref{enl}. By
Proposition~\ref{prop2} we have the estimate
\begin{equation*}
T^F_\O(\ox-tq)\le r+t\;\mbox{ whenever }\;q\in F\;\mbox{ and
}\;t\ge 0.
\end{equation*}
Since $x:=\ox-tq\in\ox+\dd\B$ when $t$ is sufficiently small, it
follows that
\begin{align*}
\la x^*,-tq\ra&\le T^F_\O(\ox-tq)-T^F_\O(\ox)+t(\ve+\gamma)\|q\|\\
&\le t+t(\ve+\gamma)\|F\|.
\end{align*}
Letting $\gamma\dn 0$ yields that $\sigma_F(-x^*)\le 1+\ve\|F\|$,
which is the upper estimate in \eqref{supp}.

To derive the lower estimate in \eqref{supp}, consider a sequence
of $\nu_k\dn 0$ as $k\to\infty$ and for any $k\in\N$ find $t_k\ge
0$ such that
\begin{equation*}
r\le t_k<r+\nu^2_k\;\mbox{ and }\;(\ox+t_kF)\cap\O\ne\emp.
\end{equation*}
The latter implies there existence of $q_k\in F$ and $w_k\in\O$
satisfying
\begin{equation*}
w_k=\ox+t_kq_k=\ox+\nu_kq_k+(t_k-\nu_k)q_k\;\mbox{ and
}\;T^F_\O(\ox+\nu_kq_n)\le t_k-\nu_k.
\end{equation*}
Moreover, we have $x_k:=\ox+t_kq_k\in\ox+\dd\B$ when $k$ is
sufficiently large. This yields
\begin{align*}
\la x^*,\nu_kq_k\ra&\le T^F_\O(\ox+\nu_kq_k)-T^F_\O(\ox)+(\ve+\gamma)\nu_k\|q_k\|\\
&\le t_k-\nu_k-r+(\ve+\gamma)\nu_k\|F\|\\
&\le\nu_k^2-\nu_k+(\ve+\gamma)\nu_k\|F\|
\end{align*}
and justifies therefore that
\begin{equation*}
1-\nu_k-(\ve+\gamma)\|F\|\le\la-x^*,q_k\ra\le\sigma_F(-x^*).
\end{equation*}
Thus we get $1-\ve\|F\|\le\sigma_F(-x^*)$ by letting $\nu_k\dn 0$
as $k\to\infty$ and taking into account that $\gg>0$ was chosen
arbitrarily. This completes the proof of the proposition.
$\h$\vspace*{0.05in}

Next we obtain an upper estimate of the $\ve$-subdifferentials of
the minimal time function \eqref{mt} at out-of-set points via the
sets of $\ve$-normals \eqref{1.2.1} to $\O$ at (generalized) {\em
projection points} and the Minkowski gauge of the dynamics
\eqref{mk}.

\begin{Proposition}\label{mink2} {\bf (upper estimate of
$\ve$-subgradients of minimal time functions at out-of-set points
via projections on targets).} Let $\ox\notin\O$ with
$T^F_\O(\ox)<\infty$, and let $\Pi^F_\O(\ox)\ne\emp$. Then for any
$\ow\in\Pi^F_\O(\ox)$ and $\ve\ge 0$ we have the estimate
\begin{eqnarray}\label{m1}
\begin{array}{ll}
\Hat\partial_\ve T^F_\O(\ox)\subset-\Hat\partial_\ve \rho_F(\bar
w-\ox)\cap \Hat N_\ve(\bar w;\O).
\end{array}
\end{eqnarray}
\end{Proposition}
{\bf Proof.} Fix a number $\ve\ge 0$ and an $\ve$-subgradient
$x^*\in\Hat\partial_\ve T^F_\O(\ox)$. Then picking any number
$\eta>0$ and employing \eqref{2.3}, we find $\delta>0$ such that
\begin{eqnarray}\label{subgradient}
\la x^*,x-\ox\ra\le
T^F_\O(x)-T^F_\O(\ox)+(\ve+\eta)\|x-\ox\|\;\mbox{ whenever
}\;x\in\ox+\delta\B.
\end{eqnarray}
Let us first show that, taking any projection point
$\ow\in\Pi^F_\O(\ox)$, we have the upper estimate
\begin{eqnarray*}
\Hat\partial_\ve T^F_\O(\ox;\O)\subset\Hat N_\ve(\bar w;\O)
\end{eqnarray*}
via $\ve$-normals \eqref{1.2.1} to the target $\O$. Indeed, fix
$\ow\in\Pi^F_\O(\ox)$ and observe that $\ow\in\O\cap(\ox+tF)$ with
$t:=T^F_\O(\ox)>0$. Hence $w\in\O\cap (w-\bar w+\ox+tF)$ for any
$w\in\O$. Specifying further $w\in\ow+\dd\B$ with $\dd>0$ from
\eqref{subgradient} and taking into account that $w-\bar w+\ox\in
\ox+\delta\B$ and $T^F_\O(w-\bar w+\ox)\le t=T^F_\O(\ox)$, we get
by \eqref{subgradient} that
\begin{align*}
\la x^*,w-\bar w\ra&\le T^F_\O(w-\bar
w+\ox)-T^F_\O(\ox)+(\ve+\eta)\|w-\bar w\|\\
&\le (\ve+\eta)\|w-\bar w\|.
\end{align*}
This implies $x^*\in\Hat N_\ve(\bar w;\O)$ by definition
\eqref{1.2.1}.

To continue the proof of estimate \eqref{m1} by involving now the
$\ve$-subdifferential of the Minkowski gauge $\rho_F$, we set
$\Tilde{x}=\bar w-\bar x$ and apply (\ref{subgradient}) with
$\ox-t(x-\tilde{x})$ and $t>0$ sufficiently small. Then
\eqref{subgradient}, Proposition~\ref{mtg1}, and the convexity of
$\rho_F$ ensure the relationships
\begin{align*}
\la x^*,-t(x-\tilde{x})\ra&\le
T^F_\O(\ox-t(x-\tilde{x}))-T^F_\O(\ox)
+(\ve+\eta)t\|x-\Tilde{x}\|\\
&\le\rho_F(\bar w-\ox+t(x-\tilde{x}))-\rho_F(\bar
w-\ox)+(\ve+\eta)t\|x-\Tilde{x}\|\\
&\le\rho_F(tx +(1-t)(\bar w-\ox))-\rho_F(\bar w-\ox)+(\ve+\eta)t\|x-\Tilde{x}\|\\
&\le t\rho_F(x)+(1-t)\rho_F(\bar w-\bar x)-\rho_F(\bar w-\ox)+(\ve+\eta)t\|x-\Tilde{x}\|\\
&=t\big(\rho_F(x)-\rho_F(\tilde{x})\big)+(\ve+\eta)t\|x-\Tilde{x}\|.
\end{align*}
Thus $-x^*\in\Hat\partial_\ve\rho_F(\bar w-\ox)$, and the proof is
complete. $\h$\vspace*{0.05in}

The last assertion of this section provides a two-sided estimate
of $\ve$-subgradients of the minimal time function \eqref{mt} at
out-of-set points $\ox\in\O$ via the set of $\ve$-normals to the
target {\em enlargements} \eqref{enl} and appropriate
characteristics of the dynamics. The results obtained extend the
ones from \cite[Theorem~4.2]{mn10} derived for the
$\ve$-subdifferential $\Hat\partial_\ve T^F_\O(\ox)$ under the
interiority assumption $0\in\i F$ and those from
\cite[Theorem~4.2]{jh} given for the Fr\'echet subdifferential
$\Hat\partial T^F_\O(\ox)$ under the calmness assumption
\eqref{calm}.

In addition to \eqref{C1}, define the {\em two-sided support set}
\begin{eqnarray}\label{C2}
S^*_\ve:=\big\{x^*\in X^*\big|\;1-\ve\|F\|\le\sigma_F(-x^*)\le
1+\ve\|F\|\},\quad\ve\ge 0,
\end{eqnarray}
which reduces to $S^*:=\{x^*\in X^*|\;\sigma_F(-x^*)=1\}$ for
$\ve=0$.

\begin{Theorem}\label{Frechet out-of-set} {\bf ($\ve$-subgradients
of minimal time functions at out-of-set points via $\ve$-normals
to target enlargements).} Let $\ox\notin\O$ with
$r:=T^F_\O(\ox)<\infty$ under our standing assumptions. Then we
have the upper estimate
\begin{equation}\label{up}
\Hat\partial_\ve T^F_\O(\ox)\subset\Hat N_\ve(\ox;\O_r)\cap
S^*_\ve\;\mbox{ for all }\;\ve\ge 0.
\end{equation}
Conversely, suppose that the minimal time function $T^F_\O$ is
calm at $\ox$ with constant $\kappa$. Then for any $x^*\in\Hat
N_\ve(\ox;\O_r)\cap S^*_\ve$ and $\ve\ge 0$, we have the inclusion
\begin{equation}\label{up1}
x^*\in\Hat\partial_{\ell\ve}T^F_\O(\ox)\;\mbox{ with
}\;\ell=\ell(x^*):=1+2\|x^*\|\cdot\|F\|+2\kappa\|F\|.
\end{equation}
\end{Theorem}
{\bf Proof.} Fix $x^*\in\Hat\partial_\ve T^F_\O(\ox)$ with $\ve\ge
0$ and observe that the inclusion $x^*\in S^*_\ve$ follows from
Proposition~\ref{norm estimate}. To justify $x^*\in\Hat
N_\ve(\ox;\O_r)$, pick $\eta>0$ and find $\dd>0$ such that
inequality \eqref{subgradient} is satisfied. Since $T^F_\O(x)\le
r=T^F_\O(\ox)$ for all $x\in\O_r$, we have
\begin{eqnarray*}
T^F_\O(x)-T^F_\O(\ox)\le 0\;\mbox{ whenever
}\;x\in\O_r\cap(\ox+\dd\B),
\end{eqnarray*}
which implies therefore that \eqref{subgradient} reduces to
\begin{eqnarray*}
\la x^*,x-\ox\ra\le(\ve+\eta)\|x-\ox\|
\end{eqnarray*}
for such $x$. Thus we get $x^*\in\Hat N_\ve(\ox;\O_r)$ by
\eqref{1.2.1} and justify the upper estimate \eqref{up}.

To prove the converse inclusion \eqref{up1} under the extra
calmness assumption, pick any $x^*\in\Hat N_\ve(\ox;\O_r)\cap
S^*_\ve$ with $\ve\ge 0$ and, applying Proposition~\ref{Frechet2}
and taking into account that $ S^*_\ve\subset C^*_\ve$ and
$\mu(x^*)\le\ell(x^*)$ for $\mu(x^*)$ in \eqref{fr} and
$\ell=\ell(x^*)$ in \eqref{up1}, we get
\begin{eqnarray}\label{fr1}
x^*\in \Hat\partial_{\ell\ve}T^F_{\O_r}(\ox)\;\mbox{ with
}\;r=T^F_\O(\ox).
\end{eqnarray}
It follows from Proposition~\ref{distance estimate} that
$T^F_\O(x)=T^F_{\O_r}(x)-r$ for any $x$ with $T^F_\O(x)<\infty$
and $T^F_\O(x)\ge r$. This yields by \eqref{fr1} that
\begin{equation}\label{fr0}
\begin{array}{ll}\disp\liminf_{x\to\ox,\;T^F_\O(x)\ge
r}\dfrac{T^F_\O(x)-T^F_\O(\ox)-
\la x^*,x-\ox\ra}{\|x-\ox\|}\\
\disp\ge\liminf_{x\to\ox,\,T^F_\O(x)\ge
r}\disp\dfrac{T^F_{\O_r}(x)-T^F_{\O_r}(\ox)-\la
x^*,x-\ox\ra}{\|x-\ox\|}\ge-\ell\ve.
\end{array}
\end{equation}
To justify \eqref{up1}, it remains to prove that
\begin{equation}\label{fr2}
\liminf_{x\to \ox,\;T^F_\O(x)<r}\dfrac{T^F_\O(x)-T^F_\O(\ox)-\la
x^*,x-\ox\ra}{\|x-\ox\|}\ge-\ell\ve.
\end{equation}

To proceed, take an arbitrary number $\gg>0$ and find $\dd>0$ such
that
\begin{align}\label{epsilon normal}
\la x^*,x-\ox\ra\le(\ve+\gamma)\|x-\ox\|\;\mbox{ whenever
}\;x\in\O_r\cap(\ox+\dd\B)
\end{align}
by $x^*\in\Hat N_\ve(\ox;\O_r)$ and $|T^F_\O(x)-T^F_\O(\ox)|\le
\kappa\|x-\ox|\|$ for all $x\in\ox+\delta\B$ by the calmness
condition. Since $\sigma_F(-x^*)\ge 1-\ve\|F\|$, there is $q\in F$
such that $\la-x^*,q\ra\ge 1-\ve\|F\|-\gamma$. Fix further a point
$x\in X$ such that $T^F_\O(x)<r$ and
\begin{equation}\label{d1}
x\in\ox+\dd_1\B\;\mbox{ with
}\;\delta_1:=\dfrac{\delta}{1+\kappa\|F\|}.
\end{equation}
Denoting $t:=T^F_\O(x)$, we take a sequence of $\nu_k\dn 0$ as
$k\to\infty$ and for any $k\in\N$ find $t_k\ge 0$, $w_k\in\O$, and
$q_k\in F$ satisfying
\begin{equation*}
t\le t_k\le t+\nu_k\;\mbox{ and }\;w_k=x+t_kq_k.
\end{equation*}
It is easy to observe that
\begin{align*}
w_k=x-(r-t_k)q+(r-t_k)q+t_kq_k\subset x-(r-t_k)q+rF
\end{align*}
when $k$ is sufficiently large. Thus for such $k$ we have
\begin{equation*}
T^F_\O(x_k)\le r\;\mbox{ with }\;x_k:=x-(r-t_k)q
\end{equation*}
and, by using $r-t=T^F_\O(\ox)-T^F_\O(x)\le\kappa\|x-\ox\|$ and
the definition of $\dd_1$ in \eqref{d1}, arrive subsequently at
the upper estimates
\begin{eqnarray}\label{d2}
\begin{array}{ll}\|x_k-\ox\|&\le\|x-\ox\|+(r-t_k)\|q\|\le
\|x-\ox\|+(r-t)\|F\|\\
&\le\|x-\ox\|+\kappa\|x-\ox\|\cdot\|F\|\le
(1+\kappa\|F\|)\delta_1<\delta,
\end{array}
\end{eqnarray}
and thus $x_k\in\ox+\dd\B$ for all $k$ sufficiently large.
Plugging now $x:=x_k$ into \eqref{epsilon normal} and employing
the middle estimate in \eqref{d2}, we get
\begin{align*}
\la x^*,x-\ox\ra-(r-t_k)\la x^*,q\ra&\le(\ve+\gamma)\|x_k-\ox\|\\
&\le (\ve+\gamma)(1+\kappa\|F\|)\|x-\ox\|
\end{align*}
for the point $x$ fixed above. The latter gives by letting
$k\to\infty$ that
\begin{align*}
\la x^*,x-\ox\ra&\le(r-t)\la x^*,q\ra+(\ve+\gamma)(1+\kappa\|F\|)\|x-\ox\|\\
&\le t-r+(\ve\|F\|+\gamma)(r-t)+(\ve+\gamma)(1+\kappa\|F\|)\|x-\ox\|\\
&\le
T^F_\O(x)-T^F_\O(\ox)+\big[\kappa(\ve\|F\|+\gamma)+(\ve+\gamma)(1+\kappa\|
F\|)\big]\|x-\ox\|,
\end{align*}
which in turn implies that
\begin{align*}
\liminf_\substack{x\to\ox}{T^F_\O(x)<r}\dfrac{T^F_\O(x)-T^F_\O(\ox)-\la x^*,
x-\ox\ra}{\|x-\ox\|}\ge-(1+2\kappa\|F\|)\ve\ge-\ell\ve,
\end{align*}
since $\gg>0$ was chosen arbitrarily. Thus we get \eqref{fr2}
and, unifying it with \eqref{fr0}, justify \eqref{up1} and
complete the proof of the theorem. $\h$

\section{Evaluating Basic and Singular Subdifferentials of Minimal
Time Functions at In-set Points of General Targets}
\setcounter{equation}{0}

In this section we obtain various formulas of inclusion and equality types
for efficient evaluations of both {\em basic} \eqref{2.4} and {\em singular}
\eqref{2.4.1} subdifferentials of minimal time functions at {\em in-set} points
$\ox\in\O$ of general nonconvex target sets $\O$. \vspace*{0.05in}

Recall that a function $\ph\colon X^*\to\oR$ is {\em sequentially weak$^*$
continuous at} $x^*$ if for any sequence $x^*_k\xrightarrow{w^*}x^*$
we have $\ph(x_k^*)\to\ph(x^*)$ as $k\to\infty$. The function $\ph$ is
sequentially  weak$^*$ continuous {\em on} a subset $S\subset X^*$ if
it has this property at each point of $S$.

In what follows we exploit the sequential weak$^*$ continuity of the dynamics
support function \eqref{sf}, which is automatic in finite dimensions due to
the following simple observation.

\begin{Proposition}\label{sup-cont} {\bf (Lipschitz continuity of support functions).}
Let $F$ be a bounded subset of a normed space $X$, and let $\sigma_F$ be the
associated support function \eqref{sf}. Then
\begin{equation}\label{supp1}
|\sigma_F(x^*_1)-\sigma_F(x^*_2)|\le\|F\|\cdot\|x^*_1-x^*_2\|\;\mbox{ for any }
\;x^*_1,x^*_2\in X^*,
\end{equation}
i.e., $\sigma_F$ is globally Lipschitz continuous with constant $\|F\|$ in the norm topology of $X^*$.
\end{Proposition}
{\bf Proof.} Fix $x^*_1,x^*_2\in X^*$ and for any $\eta>0$ find by \eqref{sf} such $q\in F$ that
$\sigma_F(x^*_1)-\eta\le\la x^*_1,q\ra$. Then we immediately have the estimates
\begin{align*}
\sigma_F(x^*_1)-\sigma_F(x^*_2)&\le\la x^*_1,q\ra-\sigma_F(x^*_2)+\eta\\
&\le\la x^*_1,q\ra-\la x^*_2,q\ra+\eta\\
&\le\|x^*_1-x^*_2\|\cdot\|F\|+\eta,
\end{align*}
which imply in turn that $\sigma_F(x^*_1)-\sigma_F(x^*_2)\le\|F\|\cdot\|x^*_1-x^*_2\|$,
since $\eta>0$ was chosen arbitrarily. Interchanging the role of $x^*_1$ and $x^*_2$ in
the latter estimate gives us \eqref{supp1}. $\h$\vspace*{0.05in}

Let us now establish two-sided relationships between the
basic subdifferential of \eqref{mt} and the basic normal to the
target in the {\em in-set} setting. The following theorem is new even for the case of
$0\in\i F$ in finite dimensions; cf. \cite[Theorem~3.6]{mn10}.

\begin{Theorem}\label{lim-min} {\bf (basic subgradients of minimal
time functions and basic normals to targets at in-set points).}
Let $\ox\in\O$ with $T^F_\O(\ox)<\infty$ for the minimal time
function \eqref{mt}, and let $C^*$ be defined in \eqref{C1} as $\ve=0$.
Then we have the upper estimate
\begin{eqnarray}\label{limiting1}
\partial T^F_\O(\ox)\subset N(\ox;\O)\cap C^*,
\end{eqnarray}
which holds as equality when the dynamics support function \eqref{sf}
is sequentially weak$^*$ continuous on the set $-[N(\ox;\O)\cap C^*]$; in particular,
when $\dim X<\infty$. If in addition $0\in F$, then we have the normal cone representation
\begin{eqnarray}\label{limiting2}
N(\ox;\O)=\bigcup_{\lambda>0}\lambda\partial T^F_\O(\ox).
\end{eqnarray}
\end{Theorem}
{\bf Proof.} To justify the upper estimate \eqref{limiting1}, fix an
arbitrary basic subgradient $x^*\in\partial T^F_\O(\ox;\O)$ and by
definition \eqref{2.4} find sequences
$\ve_k\dn 0$, $x_k\rightarrow\ox$, $T^F_\O(x_k)\to T^F_\O(\ox)=0$,
and $x^*_k\xrightarrow{w^*}x^*$ as $k\to\infty$ such that
$x^*_k\in\Hat\partial_{\ve_k} T^F_\O(x_k)$ for all $k\in\N$. If
there is a subsequence of $\{x_k\}$ (with no relabeling) that
belongs to $\O$, then we get $x^*_k\in\Hat N_{\ve_k}(x_k;\O)$ and
\begin{eqnarray}\label{m-lim}
\sigma_F(-x^*_k)\le 1+\ve_k\|F\|
\end{eqnarray}
by Proposition~\ref{Frechet1}. Passing there to the limit as
$k\to\infty$ and employing definition \eqref{2.3.1} of the basic
normal cone give us $x^*\in N(\ox;\O)$. Since furthermore
\begin{eqnarray*}
\la-x^*_k,v\ra\le 1+\ve_k\|F\|\;\mbox{ for all }\;v\in F,
\end{eqnarray*}
it follows from \eqref{m-lim} as $k\to\infty$ that
$\la-x^*,v\ra\le 1$, which justifies \eqref{limiting1} when
$\{x_k\}\subset\O$.

Consider now the other case when $x_k\notin\O$ for all $k\in\N$
sufficiently large and find by Theorem~\ref{out-of-set Frechet}
a sequence $\{w_k\}\subset\O$ satisfying
\begin{eqnarray}\label{m-lim1}
x^*_k\in\Hat N_{\ve_k+1/k}(w_k;\O)\;\mbox{ and
}\;\|x_k-w_k\|\le\|F\|T^F_\O(x_k)+1/k,\quad k\in\N.
\end{eqnarray}
Since $T^F_\O(x_k)\to 0$, it follows from the inequalities in
\eqref{m-lim1} that $w_k\to\ox$ as $k\to\infty$, and thus $x^*\in
N(\ox;\O)$ by passing to the limit in the inclusions of
\eqref{m-lim1}. We also get from Proposition~\ref{norm estimate}
that $\sigma_F(-x^*_k)\le 1+\ve_k\|F\|$ in this case, which yields
that $\sigma_F(-x^*)\le 1$ as $k\to\infty$ and completes the proof
of the upper estimate \eqref{limiting1}.\vspace*{0.05in}

Let us next justify the opposite inclusion in \eqref{limiting1}
under the additional assumption made. Pick any $x^*\in N(\ox;\O)\cap
C^*$ and by definition \eqref{2.3.1} find sequences $\ve_k\dn 0$,
$x_k\xrightarrow{\O}\ox$, and $x^*_k\xrightarrow{w^*} x^*$ such that
$x^*_k\in\Hat N_{\ve_k}(x_k;\O)$ and $\sigma_F(-x^*)\le 1$ for all
$k\in\N$. Invoking the assumed sequential weak$^*$ continuity of
$\sigma_F$ on $-[N(\ox;\O)\cap C^*]$, we get the convergence
$\sigma_F(-x^*_k)\to\sigma_F(-x^*)$ as $k\to\infty$. If
$\sigma_F(-x^*)<1$, then $\sigma_F(-x^*_k)<1$ for all large $k$.
Proposition~\ref{Frechet2} gives us a sequence $\Tilde{\ve}_k\dn 0$
such that $x^*_k\in \Hat\partial_{\tilde{\ve}_k}T^F_\O(x_k)$; hence
$x^*\in\partial T^F_\O(\ox)$.

In the other case of $\sigma_F(-x^*)=1$, denote $\gamma_k:=\sigma_F(-x^*_k)$
and get by the assumed weak$^*$ continuity that $\gamma_k\to 1$ as $k\to\infty$.
Then we have
\begin{equation}\label{m-lim0}
\dfrac{x^*_k}{\gamma}_k\in\Hat N_{\ve_k/\gamma_k}(x_k)\cap C^*\;\mbox{ and then }\;
\dfrac{x^*_k}{\gamma_k}\in\Hat\partial_{\tilde{\ve}_k}T^F_\O(x_k)
\end{equation}
for some sequence $\tilde\ve_k\dn 0$, which exists by Proposition~\ref{Frechet2}.
Passing to the limit in \eqref{m-lim0} as $k\to\infty$ yields $x^*\in\partial T^F_\O(\ox)$
and completes the proof of equality in \eqref{limiting1}.\vspace*{0.05in}

Let us finally justify representation \eqref{limiting2}.
It immediately follows from the upper estimate \eqref{limiting1}
that the inclusion ``$\subset$" holds in \eqref{limiting2}.
It remains to show that under the additional assumption $0\in F$
the opposite inclusion
\begin{eqnarray*}
N(\ox;\O)\subset\bigcup_{\lambda>0}\lambda\partial
T^F_\O(\ox),\quad\ox\in\O
\end{eqnarray*}
is satisfied. To proceed, fix any basic normal $x^*\in N(\ox;\O)$
and find by \eqref{2.3.1} sequences $\ve_k\dn 0$,
$w_k\xrightarrow{\O}\ox$, and $x^*_k\xrightarrow{w^*}x^*$ as
$k\to\infty$ such that $x_k^*\in\Hat N_{\ve_k}(w_k;\O)$ for all
$k\in\N$. Let
\begin{eqnarray*}
\lambda_k:=\sigma_F(-x^*_k)+1=\disp\sup_{v\in F}\la-x^*_k,v\ra+1,\quad k\in\N.
\end{eqnarray*}
and observe from $0\in F$ that $\lambda_k\ge 1$ for every $k$.
Moreover, the sequence $\{\lm_k\}$ is bounded in $\R$ due to the
boundedness of $F$ in $X$ and the boundedness of the weak$^*$
convergence sequence $\{x^*_k\}$ in $X^*$ by the uniform
boundedness principle. Without loss of generality, suppose that
$\lm_k\to\lm>0$ as $k\to\infty$. Then
\begin{eqnarray}\label{m-lim2}
\dfrac{x^*_k}{\lambda_k}\in\Hat\partial_{\alpha_k\ve_k/\lambda_k}T^F_\O(w_k),\quad
k\in\N,
\end{eqnarray}
with $\al_k:=2\|F\|\cdot\|x^*_k/\lambda_k\|+1\ge 1$ for all $k$.
This implies that
\begin{eqnarray*}
x^*\in\lambda\partial T^F_\O(\ox)
\end{eqnarray*}
by passing to the limit in (\ref{m-lim2}), which completes the
proof of the theorem. $\h$\vspace*{0.05in}

The next theorem provides an upper estimate of the singular
subdifferential of (non-Lipschitzian) minimal time functions at
in-set points and also justifies a case of equality therein.
As mentioned in the Introduction, the latter subdifferential
has never been considered in the literature for minimal time functions
while it is important for applications.

\begin{Theorem}\label{sin-min} {\bf (singular subgradients of minimal
time functions via basic normals to targets at in-set points).}
Define the positive dual cone of the dynamics in \eqref{mt} by
\begin{eqnarray}\label{pos}
F^*_+:=\big\{x^*\in X^*\big|\;\la x^*,v\ra\ge 0\;\mbox{ for all
}\; v\in F\big\}.
\end{eqnarray}
Then for any in-set point $\ox\in\O$ with $T^F_\O(\ox)<0$ we have
the upper estimate
\begin{equation}\label{pos1}
\partial^\infty T^F_\O(\ox)\subset N(\ox;\O)\cap F^*_+.
\end{equation}
Moreover, \eqref{pos1} holds as equality when $0\in F$ and the support
function $\sigma_F$  in \eqref{sf} is weak$^*$ continuous on the set
$-[N(\ox;\O)\cap F^*_+]$.
\end{Theorem}
{\bf Proof.} To justify \eqref{pos1}, fix any
$x^*\in\partial^\infty T^F_\O(\ox)$ and by definition \eqref{2.4.1}
find sequences $\lambda_k\dn 0$, $x_k\rightarrow\ox$, $\ve_k\dn 0$,
and $x^*_k\in\Hat\partial T^F_\O(x_k)$ such that $T^F_\O(x_k)\rightarrow
T^F_\O(\ox)=0$ and
\begin{equation*}
\lambda_kx^*_k\xrightarrow{w^*}x^*\;\mbox{ as }\;k\to\infty.
\end{equation*}
In the case of $x_k\in\O$ for a subsequence of $k\in\N$ (no
relabeling) we have
\begin{equation*}
x^*_k\in\Hat N_{\ve_k}(x_k;\O)\;\mbox{ and }\;\sigma_F(-x^*_k)\le
1+\ve_k\|F\|,\quad k\in\N,
\end{equation*}
which implies by construction \eqref{1.2.1} and
Proposition~\ref{Frechet1} that $\lambda_kx^*_k\in\Hat
N_{\lambda_k\ve_k}(x_k;\O)$ and
\begin{equation*}
\la -\lambda_kx^*_k,v\ra\le\lambda_k+\lambda_k\ve_k\|F\|\;\mbox{
whenever }\;v\in F.
\end{equation*}
By passing to the limit in the latter relationships as
$k\to\infty$, we get that $x^*\in N(\ox;\O)$ and $\la-x^*,v\ra\le 0$ for all
$v\in F$, respectively. This justifies \eqref{pos1} in the case
under consideration.

In the other case of $x_k\notin\O$ for all large $k$,
we proceed similarly to the above with using
Theorem~\ref{out-of-set Frechet} and Proposition~\ref{norm
estimate} for out-set points instead of Proposition~\ref{Frechet1}
for $x_k\in\O$; cf.\ also the proof of Theorem~\ref{lim-min}. In
this way we fully justify the upper estimate \eqref{pos1}.\vspace*{0.05in}

Let us finally prove the opposite inclusion in \eqref{pos1} under the additional
assumptions made. Fix any $x^*\in N(\ox;\O)\cap F^*_+$ and by definition \eqref{2.3.1}
find sequences $\ve\dn 0$, $x_k\xrightarrow{\O}\ox$, and $x^*_k\xrightarrow{w^*}x^*$
such that $x^*_k\in\Hat N_{\ve_k}(x_k;\O)$. We have furthermore that $\sigma(-x^*)=0$
due to $0\in F$ and $x^*\in F^*_+$. It follows from the assumed sequential weak$^*$
continuity of the support function $\sigma_F$ that $0\le\sigma_F(-x_k^*)\to
\sigma_F(-x^*)=0$. Set now
\begin{eqnarray*}
\lambda_k:=\sigma_F(-x^*_k)+\sqrt[4]{\ve_k}+1/k,\quad k\in\N,
\end{eqnarray*}
and observe that $\lambda_k\dn 0$ as $k\to\infty$ and $x^*_k/\lambda_k\in\Hat
N_{\ve_k/\lambda_k}(x_k;\O)\cap C^*$. Since $\ve_k/\lambda_k\dn 0$, by
Proposition~\ref{Frechet2} we find a sequence $\tilde{\ve}_k\dn 0$ such that
$\dfrac{x^*_k}{\lambda_k}\in\Hat\partial_{\tilde{\ve}_k}T^F_\O(x_k)$ for all $k\in\N$,
and hence $x^*\in\partial^\infty T^F_\O(\ox)$ by passing to the limit as $k\to\infty$.
This ensures the equality in \eqref{pos1} under all the assumptions made and
thus complete the proof of the theorem. $\h$\vspace*{0.05in}

Finally in this section, let us illustrate the results of Theorems~\ref{lim-min} and
\ref{sin-min}, by the following example of a two-dimensional minimal
time problem \eqref{mtp} with a nonconvex target set $\O$ and a convex
dynamics set $F$ of empty interior. In this case the minimal time function
\eqref{mt} is {\em non-Lipschitzian} and {\em nonconvex}.

\begin{Example}\label{ex-nonconvex} {\bf (basic and singular subgradients of
nonconvex and non-Lipschitzian minimal time functions at in-set points).}
{\rm Consider the convex dynamics set $F:=[-1,1]\times\{0\}
\subset\R^2$ with $\i F=\emp$ and the nonconvex target set $\O:=\R^2\setminus(-1,1)
\times(-1,1)$ in the minimal time problem \eqref{mtp}. Then the Minkowski gauge
\eqref{mk} and the minimal time function \eqref{mt} are computed, respectively, by
\begin{eqnarray*}
\rho_F(x)=\left\{\begin{array}{ll}
|x_1|&\mbox{ if }\;x\in\R\times\{0\},\\
\infty & \text{otherwise};
\end{array}\right.
\end{eqnarray*}
\begin{eqnarray}\label{mte}
T^F_{\O}(x)=\left\{\begin{array}{ll}
0&\mbox{ if }\;x\in\O,\\
1-|x_1|&\mbox{ if }\;x\notin\O.
\end{array}\right.
\end{eqnarray}
We first verify Theorem~\ref{lim-min} at the in-set point $\ox=(1,0)\in\O$. It is easy to
see that $\partial
T^F_{\O}(\ox)=[-1,0]\times\{0\}$ and that $\sigma_F(v)=|v_1|$ for any $v=(v_1,v_2)\in\R^2$.
Then
\begin{equation*}
N(\ox;\O)\cap C^*=N(\ox;\O)\cap\big\{v\in\R^2\big|\;\sigma(-v)\le 1\big\}=[-1,0]\times\{0\},
\end{equation*}
and thus \eqref{limiting1} holds as {\em equality} as well as that of \eqref{limiting2}.
We can check further the fulfillment of \eqref{pos1} as equality in Theorem~\ref{sin-min},
which yields therefore that $\partial^\infty T^F_\O(\ox)=\{0\}$. Due the result mentioned
at the end of Section~2, the latter condition fully characterizes the local Lipschitzian
property of $T^F_\O$ around $\ox$, which can be seen directly from the explicit formula
for the minimal time function given above.

Taking next another in-set point $\oy=(0,1)\in\O$, we similarly check the fulfillment of
\eqref{limiting1} and \eqref{pos1} hold as {\em equalities} with  $\partial T^F_{\O}(\oy)
=\{0\}\times\R^-$ and $\partial^\infty T^F_\O(\oy)=\{0\}\times\R^-$. The latter confirms
that $T^F_\O$ is non-Lipschitzian around $(0,1)$. We see from the precise formula \eqref{mte}
for $T^F_\O$ that this function is in fact discontinuous at $(0,1)$.}
\end{Example}

\section{Evaluating Basic and Singular Subdifferentials of Minimal
Time Functions at Out-of-set Points of General Targets}
\setcounter{equation}{0}

This section is devoted to evaluating  the basic and singular subdifferentials of the
minimal time function \eqref{mt} at {\em out-of-set} points
$\ox\notin\O$. We derive two types of results in this direction: via projection points
to the target $\O$ and via enlargements $\O_r$.

Focusing first on results of the projection type, we introduce and apply the following
property of {\em well-posedness} for minimal time functions.

\begin{Definition}\label{well} {\bf (well-posedness of minimal time
functions).} We say that the minimal time function \eqref{mt} is
{\sc well posed} at $\ox\notin\O$ with $T^F_\O(\ox)<\infty$ if for
any sequence $x_k\to\ox$ as $k\to\infty$ with
$T^F_\O(x_k)\rightarrow T^F_\O(\ox)$ there is a sequence of
projection points $w_k\in\Pi^F_\O(x_k)$ containing a convergent
subsequence.
\end{Definition}

The next proposition lists some conditions ensuring the
well-posedness of \eqref{mt}. Recall that a norm on $X$ is {\em
Kadec} if the weak and strong (with respect to this norm)
convergences agree on the boundary of the unit sphere of $X$.

\begin{Proposition}\label{well-suf} {\bf (sufficient conditions
for well-posedness).} The minimal time function \eqref{mt} is well
posed at $\ox\notin\O$ under one of the following conditions:

{\bf (a)} The target $\O$ is a compact subset of $X$;

{\bf (b)} The space $X$ is finite-dimensional and $\O$ is a closed
subset of $X$;

{\bf (c)} $X$ is reflexive, $\O\subset X$ is closed and convex,
and the Minkowski gauge \eqref{mk} generates an equivalent Kadec
norm on $X$.
\end{Proposition}
{\bf Proof.} The well-posedness of \eqref{mt} under one of the
conditions (a) and (b) is obvious. Let us justify it under
condition (c). To proceed, fix a convergent sequence $x_k\to\ox$
and observe that the property $T^F_\O(x_k)\rightarrow T^F_\O(\ox)$
is automatic when $\rho_F$ generates a norm. It is well-known in
this case that $\Pi^F_\O(x)\ne\emp$ for every $x\in X$ due to the
convexity of $\O$ and the reflexivity of $X$. Pick any $w_k\in
\Pi^F_\O(x_k)$ and observe that
\begin{equation}\label{wel}
T^F_\O(x_k)=\rho_F(x_k-w_k),\quad k\in\N.
\end{equation}
It follows that the sequence $\{w_k\}$ is bounded in $X$, and
hence---by the reflexivity of $X$---it contains a subsequence
(with no relabeling) that weakly converges to some element $\ow$.
Since $\O$ is convex and closed in $X$, it is also weakly closed;
this $\ow\in\O$. By the lower semicontinuity of $\rho_F$ in the
weak topology of $X$ and by \eqref{wel} we have the relationships
\begin{equation*}
\rho_F(\ox-\bar w)\le\liminf_{k\to\infty}\rho_F(x_k-w_k)=\liminf_{k\to\infty}
T^F_\O(w_k-x_k)=T^F_\O(\ox),
\end{equation*}
which imply that $\bar w\in\Pi^F_\O(\ox)$ and $T^F_\O(\ox-\bar
w)=\rho_F(\bar x-\bar w)$. Since $\rho_F$ generates a Kadec norm
on $X$, it follows from $\rho_F(x_k-w_k)\to\rho_F(\bar x-\bar w)$
and the weak convergence of $x_k-w_k$ to $\ox-\ow$ that in fact
the sequence $x_k-w_k$ converges strongly in $X$, and hence
$w_k\to\ow$ as $k\to\infty$. This completes the proof of the
proposition. $\h$\vspace*{0.05in}

Now we use the well-posedness property of $T^F_\O$ to derive upper
estimates of both basic and singular subdifferentials of the
minimal time function at out-of-set points.

\begin{Theorem}\label{mink3} {\bf (basic and singular subgradients
of minimal time functions at out-of-set points via projections).} Let
$\ox\notin\O$ with $T^F_\O(\ox)<\infty$, and let the minimal time
function \eqref{mt} be well posed at $\ox$. Then we have the
estimates
\begin{eqnarray}\label{wel1}
\partial T^F_\O(\ox)\disp\subset\bigcup_{\bar w\in\Pi^F_\O(\ox)}\big
[-\partial\rho_F(\bar w-\ox)\cap N(\bar w;\O)\big],
\end{eqnarray}
\begin{eqnarray}\label{wel2}
\begin{array}{ll}\disp\partial^\infty T^F_\O(\ox)&\subset\disp\bigcup_{\bar w\in
\Pi^F_\O(\ox)}\big[-\partial^\infty\rho_F(\bar w-\ox)\cap N(\bar w; \O)\big]\\\\
&\subset\disp\bigcup_{\bar w\in\Pi^F_\O(\ox)}\big[N(\bar w;\O)\cap
F^*_+\big]
\end{array}
\end{eqnarray}
with the positive dual cone $F^*_+$ of the dynamics defined in
\eqref{pos}.
\end{Theorem}
{\bf Proof.} Pick any basic subgradient $x^*\in\partial
T^F_\O(\ox)$ and by definition \eqref{2.4} find sequences
$\ve_k\dn 0$, $x_k\xrightarrow{T^F_\O}\ox$, and
$x^*_k\in\Hat\partial_{\ve_k}T^F_\O(x_k;\O)$ as
$k\rightarrow\infty$ such that $x^*_k\xrightarrow{w^*}x^*$ and
\begin{eqnarray}\label{5.2}
x^*_k\in\Hat\partial_{\ve_k}T^F_\O(x_k;\O)\;\mbox{ for all
}\;k\in\N.
\end{eqnarray}
By the well-posedness property of \eqref{mt} there is a sequence
$w_k\in\Pi^F_\O(x_k;\O)$, which contains a subsequence (no
relabeling) converging to some $\bar w$. It follows from
definitions \eqref{projection} of the generalized projection, the
convergence $T^F_\O(x_k)\to T^F_\O(\ox)$, and the assumptions made
that $\ow\in\Pi^F_\O(\ox)$. Applying Proposition~\ref{mink2} to
\eqref{5.2}, we have
\begin{eqnarray*}
x^*_k\in-\Hat\partial_{\ve_k}\rho_F(x_k-w_k)\cap\Hat
N_{\ve_k}(w_k;\O),\quad k\in\N,
\end{eqnarray*}
which yields in turn the upper estimates \eqref{wel1} by passing
to the limit as $k\to\infty$.

Let us now prove both inclusions in \eqref{wel2}. Taking an
arbitrary singular subgradient $x^*\in\partial^\infty
T^F_\O(\ox)$, find by \eqref{2.4.1} sequences $\ve_k\dn 0$,
$\lambda_k\dn 0$, $x_k\xrightarrow{T^F_\O}\ox$, and
$x^*_k\in\Hat\partial_{\ve_k}T^F_\O(x_k;\O)$ such that
$\lambda_kx^*_k\xrightarrow{w^*}x^*$ as $k\rightarrow\infty$ and
\begin{eqnarray}\label{5.2a}
x^*_k\in\Hat\partial_{\ve_k}T^F_\O(x_k)\;\mbox{ for all
}\;k\in\N.
\end{eqnarray}
By the well-posedness property of \eqref{mt} there is a sequence
$w_k\in\Pi^F_\O(x_k;\O)$ that contains a subsequence (no
relabeling) converging to some $\bar w$. As discussed above, we
have $\bar w\in\Pi^F_\O(\ox)$. Applying Proposition~\ref{mink2} to
\eqref{5.2a} allows us to conclude that
\begin{eqnarray}\label{wel3}
-\lm_kx^*_k\in\lm_k\Hat\partial_{\ve_k}\rho_F(x_k-w_k)\;\mbox{
and}\;x^*_k\in\lm_k\Hat N_{\ve_k}(w_k;\O),\quad k\in\N.
\end{eqnarray}
Letting $k\to\infty$ in both inclusions of \eqref{wel3}, we arrive
at the first estimate in \eqref{wel2}.

To justify the remaining inclusion in \eqref{wel2}, observe by the
arguments similar to the corresponding ones in
Theorem~\ref{sin-min} (cf.\ also the proof of Theorem~\ref{convex
case 2} below for more details in the like setting) that we have
the implication
\begin{eqnarray*}
-x^*_k\in\Hat\partial_{\ve_k}\rho_F(x_k-w_k)\Longrightarrow
\sigma_F(-x^*_k)\le 1+\ve_k\|F\|,\quad k\in\N.
\end{eqnarray*}
It yields by \eqref{wel3} that $x^*\in N(\ow;\O)\cap F^*_+$ similarly to
the proof of Theorem~\ref{sin-min}, which thus completes the proof of
this theorem. $\h$\vspace*{0.05in}

The following example illustrates some features of the results obtained
in Theorem~\ref{mink3}.

\begin{Example}\label{ex-nonconvex1} {\bf (basic and singular subgradients of
nonconvex and non-Lipschitzian minimal time functions at out-of-set points).}
{\rm  Consider the setting of Example~\ref{ex-nonconvex}, where the minimal time
function is computed by formula \eqref{mte}. Take the out-of-set point $\bar z=(1/2,1/2)
\notin\O$ and verify the conclusions of Theorem~\ref{mink3}. The well-posedness property
\eqref{wel} holds by Proposition~\ref{well-suf}(ii). It is easy to check that
$\Pi^F_\O(\oz)=\{\bar w\}$ with $\bar w=(1,1/2)$ for the {\em Euclidean norm}
in the projection operator \eqref{projection}. Thus we arrive at the {\em equality}
\begin{equation*}
\partial T^F_\O(\bar z)=\big\{(-1,0)\big\}=-\partial\rho_F(\bar w-\bar z)\cap N(\bar w;\O)
\end{equation*}
in \eqref{wel1} and similarly get the equality in \eqref{wel2} with $\partial^\infty
T^F_\O(\oz)=\{0\}$,
which is in accordance with the local Lipschitz property of $T^F_\O$ around this point that
obviously follows from the explicit formula \eqref{mte}. Note that both inclusions \eqref{wel1}
and \eqref{wel2} are {\em strict} in this example if the projection in \eqref{projection}
is taken with respect to the {\em maximum norm} on the plane.}
\end{Example}

Let us further address a natural question about getting counterparts of Theorems~\ref{lim-min}
and \ref{sin-min} on upper estimates for basic and singular subgradients of
the minimal time function \eqref{mt} at {\em out-of-set} points via basic normals
to the enlargement $\O_r$ of the target set $\O$. However, simple examples
show the failure of such estimates. For instance, consider the minimal time problem
\eqref{mtp} in $X=\R^2$ with $F=\B$ and
$\O=\{x\in\R^2|\;\|x\|\ge 1\}$. Then for $\ox=0$ and
$r=T^F_\O(\ox)=1$ we have $N(\ox;\O_r)=\{0\}$ while $\partial
T^F_\O(\ox)=\{x\in\R^2|\;\|x\|=1\}$.

It occurs that the appropriate analogs of the upper estimates in Theorem~\ref{lim-min}
and \ref{sin-min} hold at $\ox\notin\O$ with the replacement of $\partial T^F_\O(\ox)$
and $\partial^\infty T^F_\O(\ox)$ therein by the {\em one-sided} modifications of these
constructions for $\ph=T^F_\O$ defined by
\begin{eqnarray}\label{bas-os}
\partial_{\ge}\ph(\ox):=\disp\Limsup_\substack{x\st{\ph^+}{\to}\ox}{\ve\dn 0}\Hat\partial_\ve\ph(x),
\end{eqnarray}
\begin{eqnarray}\label{sin-os}
\partial_{\ge}^\infty\ph(\ox):=\disp\Limsup_\subsubstack{x\st{\ph^+}{\to}\ox}{\ve\dn 0}{\lm\dn 0}
\lm\Hat\partial_\ve\ph(x),
\end{eqnarray}
where the symbol $x\st{\ph^+}{\to}\ox$ signifies that $x\to\ox$ with $\ph(x)\to\ph(\ox)$ and
$\ph(x)\ge\ph(\ox)$. Note that the basic one-sided construction \eqref{bas-os} was introduced in
\cite{bmnam} and applied therein to the study of distance function (see also \cite[Sec.1.3.3]{mor06a}
and \cite{mn10}) while the singular one \eqref{sin-os} appears here for the first time. Observe that
we always have the inclusions
\begin{equation*}
\Hat\partial\ph(\ox)\subset\partial_\ge\ph(\ox)\subset\partial\ph(\ox)\;\mbox{ and }\;\partial^
\infty_\ge\ph(\ox)\subset\partial^\infty\ph(\ox)
\end{equation*}
which show, in particular, that $\partial_{\ge}\ph(\ox)=\partial\ph(\ox)$ if $\ph$ is
{\em subdifferentially regular} at $\ox$, i.e., $\Hat\partial\ph(\ox)=\partial\ph(\ox)$; the latter is
always the case for convex function.

Now we are ready to establish the corresponding counterparts of Theorem~\ref{lim-min} and \ref{sin-min}
at out-of-set points by using the one-sided constructions \eqref{bas-os} and \eqref{sin-os}.

\begin{Theorem}\label{one-sided} {\bf (one-sided basic and singular subgradients of minimal time
functions at out-of-set points).} Let the minimal time function $T^F_\O$ be continuous around some
point $\ox\notin\O$, let $r=T^F_\O(\ox)$, and let the sets $C^*$, $S^*$, and $F^*_+$ be defined in \eqref{C1},
\eqref{C2}, and \eqref{pos}, respectively. Then we have the upper estimates
\begin{equation}\label{oc}
\partial_\ge T^F_\O(\ox)\subset N(\ox;\O_r)\cap C^*\;\mbox{ and }\;
\partial^\infty_\ge T^F_\O(\ox)\subset N(\ox;\O_r)\cap F^*_+,
\end{equation}
where the first one can be replaced by the equality
\begin{eqnarray}\label{oc0}
\partial_\ge T^F_\O(\ox)=N(\ox;\O_r)\cap S^*
\end{eqnarray}
if the support function $\sigma_F$ is sequentially weak$^*$
continuous on the set $-[N(\ox;\O_r)\cap C^*]$ and if $T^F_\O$ is locally Lipschitzian around $\ox$.
Furthermore, the normal cone representation
\begin{equation}\label{equality}
N(\ox;\O_r)=\bigcup_{\lambda\ge 0}\lambda\partial_\ge T^F_\O(\ox)
\end{equation}
takes place with the convention $0\times\emp=0$ provided that $0\in\i F$.
\end{Theorem}
{\bf Proof.} We justify only the first inclusion in \eqref{oc}; the second one is proved similarly by
taking into account the proof of Theorem~\ref{sin-min}. To proceed, pick any $x^*\in
\partial_\ge T^F_\O(\ox)$ and by \eqref{bas-os} find sequences $\ve_k\dn 0$,
$x_k\st{\ph^+}{\to}\ox$, and $x^*_k\st{w^*}{\to}x^*$ as $k\to\infty$ such that
\begin{eqnarray*}
x^*_k\in\Hat\partial_{\ve_k}T^F_\O(x_k)\;\mbox{ for all }\;k\in\N.
\end{eqnarray*}
If $T^F_\O(x_k)=r$ for some subsequence of $\{k\}$, we have by the upper estimate \eqref{up} of
Theorem~\ref{Frechet out-of-set} the relationships
\begin{eqnarray*}
x^*_k\in\Hat N_{\ve_k}(x_k;\O_r)\;\mbox{ and }\;1-\ve_k\|F\|\le\sigma_F(-x^*_k)\le 1+\ve_k\|F\|
\end{eqnarray*}
held along this subsequence. Passing there to the limit as $k\to\infty$ gives us the inclusions
$x^*\in N(\ox;\O_r)$ and $x^*\in C^*$, which justify the first estimate in \eqref{oc} in this case
even without the continuity assumption on the minimal time function.

In the other case of $T^F_\O(x_k)>r$ for all $k\in\N$ sufficiently large, the assumed
continuity of $T^F_\O$ ensures that for such $k$ we have that $T^F_\O(x)>r$ whenever $x$ is near $x_k$.
Employing then Proposition~\ref{distance estimate} ensures the equality
\begin{eqnarray*}
T^F_\O(x)=r+T^F_{\O_r}(x)\;\mbox{ for all $x$ near $x_k$}.
\end{eqnarray*}
The latter implies by definition \eqref{2.3} that
\begin{eqnarray*}
x^*_k\in\Hat\partial_{\ve_k}T^F_\O(x_k)=\Hat\partial_{\ve_k}T^F_{\O_r}(x_k),\quad k\in\N.
\end{eqnarray*}
The rest of the proof of the first inclusion in \eqref{oc} follows the arguments in the proof of
Theorem~\ref{lim-min}, which in turn are based on the variational result of Theorem~\ref{out-of-set Frechet}.

Let us next justify equality \eqref{oc0} provided the fulfillment of the additional
weak$^*$ continuity and Lipschitzian assumptions made in the theorem. It follows from the proof above that
the latter assumption implies the inclusion ``$\subset$" in \eqref{oc0}. To justify the opposite inclusion
``$\supset$" therein, fix any $x^*\in N(\ox;\O_r)\cap S^*$ and find by \eqref{2.3.1} sequences $\ve_k\dn 0$,
$x_k\xrightarrow{\O_r}\ox$, and $x^*_k\xrightarrow{w^*}x^*$ as
$k\to\infty$ with $x^*_k\in\Hat N_{\ve_k}(x_k;\O_r)$, $k\in\N$. The sequential weak$^*$ continuity of
$\sigma_F$ at $-x^*$ ensures that
\begin{eqnarray*}
\gamma_k:=\sigma_F(-x^*_k)\to\sigma_F(-x^*)=1\;\mbox{ as }\;k\to\infty.
\end{eqnarray*}
By the definition of $S^*$ in \eqref{C2} we may assume with no lost of generality that
\begin{equation}\label{oc1}
\dfrac{x^*_k}{\gamma_k}\in\Hat N_{\ve_k/\gamma_k}(x_k;\O_r)\cap S^*\;\mbox{ for all }\;k\in\N.
\end{equation}
It follows further that $T^F_\O(x_k)=r$ for large $k$, since the opposite assumption on $T^F_\O(x_k)<r$
implies by the continuity of $T^F_\O$ that $x_k\in\i\O_r$, which contradicts the condition $x^*\ne 0$
held by \eqref{oc1}. Employing the second part of Theorem~\ref{Frechet out-of-set}, find a sequence
$\tilde{\ve}_k\dn 0$ such that
\begin{eqnarray*}
\dfrac{x^*_k}{\gamma_k}\in\Hat\partial_{\tilde{\ve}_k}T^F_\O(x_k)\;\mbox{ for all }\;k\in\N.
\end{eqnarray*}
Passing there to the limit as $k\to\infty$ and using definition \eqref{bas-os} justify equality \eqref{oc0}.

Let us finally prove representation \eqref{equality} correcting the corresponding arguments given in
\cite[Theorem~4.4]{mn10}. Note that the inclusion ``$\subset$" in \eqref{equality} follows from the first
inclusion \eqref{oc} and the cone property  of $N(\ox;\O_r)$. To prove the opposite inclusion $\supset$" in
\eqref{equality}, fix any $x^*\in N(\ox;\O_r)$ and assume that $x^*\ne 0$, since otherwise $x^*$ belongs to
the right-hand side of (\ref{equality}) by our convention. In this case $\gg:=\sigma_F(-x^*)>0$ due to $0\in\i F$.
By definition \eqref{2.3.1} of the basic normal cone, there are sequences $\ve_k\dn 0$, $x_k\xrightarrow{\O_r}\ox$,
and $x^*_k\xrightarrow{w^*}x^*$ with $x^*_k\in\Hat N_{\ve_k}(x_k;\O_r)$. By $0\in\i F$ the minimal time function
(\ref{mt}) is Lipschitz continuous and hence $T^F_\O(x_k)=r$ when $k$ is sufficiently large. Indeed, if
$T^F_\O(x_k)<r$ for a subsequence (without relabeling), then $x_k\in\i\O_r$, which implies that
$\|x^*_k\|\le\ve_k$ and leads to a contradiction by $\|x^*\|\le\liminf\|x^*_k\|$ as $k\to\infty$. Define further
$\lambda_k:=\sigma_F(-x^*_k)$ and observe by $x^*_k\xrightarrow{w^*}x^*$ that $\lambda_k\ge\gg/2>0$ for all
$k$ sufficiently large. Moreover, $\lambda_k$ is bounded, and hence we may assume that $\lambda_k\to\lambda\ge\gg/2$
as $k\to\infty$. Then
\begin{equation*}
\Tilde{x}^*_k:=\dfrac{x^*_k}{\lambda_k}\in\Hat N_{\ve_k/\lambda_k}(x_k)\;\mbox{ and }\;\sigma_F(-\Tilde{x}^*_k)=1,
\end{equation*}
which yields by Theorem~\ref{Frechet out-of-set} that $\Tilde{x}^*_k\in\Hat\partial_{\tilde{\ve}_k}T^F_{\O}(x_k)$ with
$\tilde{\ve}_k\to 0$ as $k\to \infty$. The latter implies the inclusions
\begin{eqnarray*}
x^*\in \lambda\partial_\ge\lambda T^F_\O(\ox)\subset\bigcup_{\lambda\ge 0}\lambda\partial_\ge T^F_\O(\ox),
\end{eqnarray*}
which justify \eqref{equality} complete the proof of the theorem. $\h$

\section{Computing Basic and Singular Subdifferentials of Convex Minimal Time Functions}
\setcounter{equation}{0}

The concluding section of the paper concerns the minimal time problem \eqref{mtp} with convex data, i.e.,
under the assumption that the target set $\O$ is a convex subset of an arbitrary Banach space $X$. By
Proposition~\ref{conv} this property is equivalent to the convexity of the minimal time function \eqref{mt}.
In what follows we add the convexity of \eqref{mt} to our standing assumptions formulated in Section~1 and
refer to this setting as to the {\em convex minimal time problem} and/or the {\em convex minimal time function}.

Due to the representations of $\ve$-normals to convex sets \eqref{n-convex} and $\ve$-subgradients
of convex functions \eqref{2.3a} we have specifications of the results obtained in Section~4 in the
case of convex minimal time functions. The same can be said regarding the results of
Sections~5 and 6 concerning the basic subdifferential and normal cone for convex functions and sets, which
reduce to those in convex analysis. We can also specify to the case of convex minimal time functions the
results derived above for the singular subdifferential; see \cite[Proposition~8.12]{rw} for its various
representations in the general framework of convex analysis.\vspace*{0.05in}

In this section we show that, besides the aforementioned specifications, the convex case allows us to obtain
{\em equalities} in the upper estimates of Sections~5 and 6 for the basic and singular subdifferentials of
\eqref{mt} at both in-set and out-of-set points with no additional assumptions in general Banach spaces. Let
us start with computing the basic subdifferential \eqref{2.3}; cf.\ Theorem~\ref{lim-min} and
Theorem~\ref{one-sided}, where $\partial_{\ge}T^F_\O(\ox)=\partial T^F_\O(\ox)$ in the convex case.

\begin{Theorem}\label{convex case} {\bf (basic subgradients of convex minimal time functions).}
Let the function $T^F_\O$ in \eqref{mt} be convex. Then the following assertions hold:

{\bf (i)} For any $\ox\in\O$ we have the representation
\begin{equation}\label{c1}
\partial T^F_\O(\ox)=N(\ox;\O)\cap C^*,
\end{equation}
where $C^*$ is defined in \eqref{C1}.

{\bf (ii)} For any $\ox\notin\O$ with $T^F_\O(\ox)<\infty$ we have the representation
\begin{equation}\label{c2}
\partial T^F_\O(\ox)=N(\ox;\O_r)\cap S^*,
\end{equation}
where $r=T^F_\O(\ox)>0$ and $S^*$ is defined in \eqref{C2}.
\end{Theorem}
{\bf Proof.} Equality \eqref{c1} in (i) follows directly from Propositions~\ref{Frechet1}
and \ref{Frechet2} with $\ve=0$ therein and the fact that $\Hat\partial T^F_\O(\ox)
=\partial T^F_\O(\ox)$ for convex functions.

To justify representation \eqref{c2} in the out-of set case (ii), observe first that the inclusion
``$\subset$" follows from the first part of Theorem~\ref{Frechet out-of-set}. It remains to prove
the converse inclusion ``$\supset$". Fix $x^*\in N(\ox;\O_r)$ with
$\sigma_F(-x^*)=1$ and show that
\begin{equation}\label{convexin}
\la x^*,x-\ox\ra\le T^F_\O(x)-T^F_\O(\ox)\;\mbox{ for all }\;x\in X.
\end{equation}
Indeed, we get from $x^*\in N(\ox;\O_r)$ and the normal cone construction for convex sets that
\begin{equation*}
\la x^*,x-\ox\ra\le 0\;\mbox{ whenever }\;x\in\O_r.
\end{equation*}
It follows from \eqref{c1} that $x^*\in\partial T^F_{\O_r}(\ox)$ and hence
\begin{equation*}
\la x^*, x-\ox\ra \le T^F_{\O_r}(x)\;\mbox{ for any }\;x\in X.
\end{equation*}
It is clear that \eqref{convexin} holds when $x\notin\O_r$, since in this
case $T^F_{\O_r}(x)=T^F_\O(x)-r$ by Proposition~\ref{distance estimate}. In the other case of
$t=T^F_\O(x)\le r$, for any $\ve>0$ sufficiently small pick $q\in F$ with $\la x^*,
-q\ra\ge 1-\ve$ and get $T^F_\O(x-(r-t)q)\le r$ by Proposition~\ref{prop2}. This gives
\begin{equation*}
\la x^*,x-\ox\ra\le(r-t)\la x^*,q\ra\le(t-r)(1-\ve),
\end{equation*}
By the arbitrary choice of $\ve>0$ the latter justifies \eqref{convexin} in this case. Thus we arrive at
$x^*\in\partial T^F_{\O_r}(\ox)$ and complete the proof of theorem. $\h$\vspace*{0.05in}

The next result provides precise representations for the singular subdifferential of the convex minimal time
function \eqref{mt} in both in-set and out-of-set cases; cf.\ Theorems~\ref{sin-min} and Theorem~\ref{one-sided},
where $\partial^\infty_\ge T^F_\O(\ox)=\partial^\infty T^F_\O(\ox)$ in the convex case.

\begin{Theorem}\label{convexsingular1}{\bf (singular subgradients of convex minimal time functions).} Let the
function $T^F_\O$ in \eqref{mt} be convex and lower semicontinuous around $\ox$, and let $F^*_+$ be defined in
\eqref{pos}. The following assertions hold:

{\bf (i)} If $\ox\in\O$, then we have
\begin{equation}\label{c3}
\partial^\infty T^F_\O(\ox)=N(\ox;\O)\cap F^*_+.
\end{equation}

{\bf (ii)} If $\ox\notin\O$ and $T^F_\O(\ox)<\infty$, then
\begin{equation}\label{c4}
\partial^\infty T^F_\O(\ox)=N(\ox;\O_r)\cap F^*_+\;\mbox{ with }\;r=T^F_\O(\ox).
\end{equation}
\end{Theorem}
{\bf Proof.} Taking into account that the subdifferential of convex analysis agrees with the Fr\'echet subdifferential
for convex functions and following the proof of \cite[Lemma~2.37]{mor06a} with replacing the fuzzy sum rule for Fr\'echet
subgradients of l.s.c.\ functions in Asplund spaces by the exact sum rule (Moreau-Rockafellar theorem) in convex analysis
in Banach spaces, we get the singular subdifferential representations under the assumptions made:
\begin{equation}\label{convex singular representation}
\partial^\infty T^F_\O(\ox)=\Limsup_\substack{x\xrightarrow{T^F_\O}\ox}{\lm\dn 0}\lambda\partial T^F_\O(x)=\big\{ x^*\in X^*\big|\;
(x^*,0)\in N\big((\ox,T^F_\O(\ox));\epi T^F_\O\big)\big\}.
\end{equation}
It is easy to check that
\begin{equation*}
\big\{ x^*\in X^*\big|\;(x^*,0)\in N\big((\ox,T^F_\O(\ox));\epi T\big)\big\}=N\big(\ox;\dom T^F_\O\big),
\end{equation*}
and hence we have by the second representation in \eqref{convex singular representation} and Theorem~\ref{sin-min} that
\begin{equation}\label{c5}
\partial^\infty T^F_\O(\ox)=N\big(\ox;\dom T^F_\O\big)\subset N(\ox;\O)\cap F^*_+.
\end{equation}
Let us now justify the opposite inclusion in \eqref{c5}, i.e.,
\begin{equation}\label{c6}
N(\ox;\O)\cap F^*_+\subset N\big(\ox;\dom T^F_\O\big).
\end{equation}
To proceed, pick arbitrary $x^*\in N(\ox;\O)\cap F^*_+$ and $x\in\dom T^F_\O$ and then find by \eqref{mtp} a number $t\ge 0$ such that
$(x+tF)\cap\O\ne\emp$. Fix further $q\in F$ and $w\in\O$ with $x+tq=w$ and obtain the relationships
\begin{align*}
\la x^*,x-\ox\ra&=\la x^*,w-tq-\ox\ra\\
&=\la x^*,w-\ox\ra-t\la x^*,q\ra\le 0,
\end{align*}
since $\la x^*,w-\ox\ra\le 0$ by $x^*\in N(\ox;\O)$ and $\la x^*,q\ra\ge 0$ by $x^*\in F^*_+$. Thus we get
\eqref{c6} and arrive at the singular subdifferential representation \eqref{c3} in the in-set case.

To justify further representation \eqref{c4} in the out-of-set case $\ox\in\O_r$ with $r=T^F_\O(\ox)$,
observe from the equality in \eqref{c5} that
\begin{equation*}
\partial^\infty T^F_\O(\ox)=N\big(\ox;\dom T^F_\O\big)\subset N(\ox;\O_r)
\end{equation*}
due to the obvious inclusions $\O_r\subset\dom T^F_\O$ and $N(\ox;\Theta_2)\subset N(\ox;\Theta_1)$
for any convex sets $\ox\in\Theta_1\subset\Theta_2$. Fix now $x^*\in \partial^\infty T^F_\O(\ox)$
and find by the first representation in \eqref{convex singular representation} sequences $x_k\xrightarrow
{T^F_\O}\ox$, $x^*_k\in\partial T^F_\O(x_k)$, and
$\lambda_k\dn 0$ such that
\begin{equation*}
\lambda_k x^*_k\xrightarrow{w^*}x^*\;\mbox{ as }\;k\to\infty.
\end{equation*}
It follows from Theorem~\ref{convex case}(ii) that $\sigma_F(-x^*_k)=1$ whenever $k\in\N$ is sufficiently
large. Hence picking any $q\in F$, we have $\la-\lambda_ kx^*_k,q\ra\le\lambda_k$ for all such $k$.
This yields $\la x^*,q\ra\ge 0$ by letting $k\to\infty$. Thus it gives $x^*\in F^*_+$ justifying the
inclusion
\begin{equation*}
\partial^\infty T^F_\O(\ox)\subset N(\ox;\O_r)\cap F^*_+.
\end{equation*}
To get \eqref{c4}, it remains to prove the converse inclusion
\begin{equation*}
N(\ox; \O_r)\cap F^+\subset N\big(\ox;\dom T^F_\O\big).
\end{equation*}
Fix $x^*\in N(\ox;\O_r)\cap F^*_+$ and pick any $x\in\dom T^F_\O$, which ensures the existence of
$t\ge 0$ such that $(x+tF)\cap\O\ne\emp$. Take $q\in F$ and $w\in\O$ satisfying $x+tq=w$. Then
\begin{align*}
\la x^*,x-\ox\ra&=\la x^*,w-tq-\ox\ra\\
&=\la x^*,w-\ox\ra-t\la x^*,q\ra\le 0
\end{align*}
by $w\in\O\subset\O_r$ and $\ox\in\O_r$, which completes the proof of the theorem. $\h$\vspace*{0.05in}

The last result of this section establishes representations of the convex subdifferential of $T^F_\O$ via
that of the Minkowski gauge; in particular, it justifies the equality in the upper estimate of $\partial
T^F_\O(\ox)$ from Theorem~\ref{mink3} at out-of-set points. Note that even the upper estimate \eqref{wel1}
itself is new with no well-posedness assumption in general Banach spaces.

\begin{Theorem}\label{convex case 2} {\bf (precise relationships between convex subdifferentials of minimal
time and Minkowski functions in out-of-set points).} Let the function $T^F_\O$ in \eqref{mt} be convex, and let
$\ox\notin\O$ be such that $\Pi^F_\O(\ox)\ne\emp$ with $r=T^F_\O(\ox)<\infty$. Then for any $\bar w\in\Pi^F_\O
(\ox)$ we have the relationships
\begin{eqnarray}\label{c7}
\begin{array}{ll}\partial T^F_\O(\ox)&=N(\ox;\O_r)\cap\big[-\partial\rho_F(\bar w-\ox)\big]\\
&\subset N(\bar w;\O)\cap\big[-\partial\rho_F(\bar w-\ox)\big].
\end{array}
\end{eqnarray}
If in addition $0\in F$, then the inclusion in \eqref{c7} holds as equality and thus
\begin{align*}
\partial T^F_\O(\ox)=N(\bar w;\O)\cap\big[-\partial\rho_F(\bar w-\ox)\big].
\end{align*}
\end{Theorem}
{\bf Proof.} It follows from Theorem~\ref{convex case}(ii) that
$\partial T^F_\O(\ox)\subset N(\ox;\O_r)$. Furthermore
\begin{equation*}
\partial T^F_\O(\ox)\subset-\partial\rho_F(\ox-\bar w)
\end{equation*}
by Proposition~\ref{mink2} as $\ve=0$, and thus
\begin{equation}\label{c8}
\partial T^F_\O(\ox)\subset N(\ox;\O_r)\cap\big[-\partial
\rho_F(\bar w-\ox)\big].
\end{equation}
To prove the opposite inclusion ``$\supset$" to \eqref{c8}, fix any $x^*\in N(\ox;\O_r)\cap
\big[-\partial\rho_F(\bar w-\ox)\big]$. By Theorem~\ref{convex case}(ii) it suffices to show that
\begin{eqnarray}\label{c9}
x^*\in S^*,\;\mbox{ i.e., }\;\sigma_F(-x^*)=1.
\end{eqnarray}
To this end, observe that $T^F_{\{0\}}(x)=\rho_F(-x)$, which implies that
\begin{eqnarray*}
-\partial\rho_F(x)=\partial T^F_{\{0\}}(-x)\;\mbox{ and hence }\;-\partial\rho_F(\bar w-\ox)=\partial T^F_{\{0\}}(\bar
x-\bar w),\quad x\in X.
\end{eqnarray*}
Since $\bar x-\bar w\notin\{0\}$, we get \eqref{c9} from Theorem~\ref{convex case}(ii) and thus justify the equality in \eqref{c7}.

Further, it is not hard to check that $\partial T^F_\O(\ox)\subset N(\bar w;\O)$ and hence
\begin{align*}
\partial T^F_\O(\ox)\subset N(\bar w;\O)\cap\big[-\partial\rho_F(\bar w-\ox)\big],
\end{align*}
which implies the inclusion in \eqref{c7}.

To finish the proof, it remains to show that
\begin{align}\label{c10}
N(\bar w;\O)\cap\big[-\partial\rho_F(\bar w-\ox)\big]\subset
N(\ox;\O_r)\cap\big[-\partial\rho_F(\bar w-\ox)\big]
\end{align}
under the additional assumption that $0\in F$ in which case we have $\rho(0)=0$. It suffices to verify
that for each $x^*\in N(\bar w;\O)\cap \big[-\partial\rho_F(\bar w-\ox)\big]$ we have $x^*\in N(\ox;\O_r)$.

To proceed, pick any $x\in \O_r$ and for an arbitrary small $\ve>0$ find $t<r+\ve$, $q\in F$,
and $w\in\O$ with $w=x+tq$. Then $\la-x^*,q\ra\le\sigma_F(-x^*)\le 1$ and
\begin{align*}
\la x^*,x-\ox\ra&=\la x^*,w-tq-\ox\ra\\
&=t\la-x^*,q\ra+\la x^*,w-\bar w\ra +\la x^*,
\bar w-\ox\ra\\
&\le t+\la x^*, w-\bar w\ra+\la x^*,\bar w-\ox\ra\\
&\le T^F_\O(\ox)+\ve+\la x^*,w-\bar w\ra+\la x^*,
\bar w-\ox\ra.
\end{align*}
We have $\la x^*,w-\bar w\ra\le 0$ due to $x^*\in N(\bar w;\O)$ and
\begin{align*}
\la x^*,\bar w-\ox\ra=\la-x^*,0-(\bar w-\ox)\ra\le
\rho_F(0)-\rho_F(\bar w-\ox)=-T^F_\O(\ox)
\end{align*}
by $-x^*\in \partial\rho_F(\bar w-\ox)$. It follows therefore that $\la x^*,x-\ox\ra\le\ve$ for all $x\in\O_r$,
and hence $x^*\in N(\ox;\O_r)$ because $\ve>0$ was chosen arbitrary small. Thus we arrive at \eqref{c10} and
complete the proof of the theorem. $\h$\vspace*{0.05in}

Finally, let us present an example that illustrates computing the basic and singular subdifferentials of convex
minimal time functions at in-set and out-of-set points.

\begin{Example}\label{convex-ex} {\bf (subdifferentiation of convex minimal time functions).} {\rm In $\R^2$, consider
the convex dynamics $F=[-1,1]\times\{0\}$ of empty interior and the convex target $\O=[-1,1]\times[-1,1]$. In this case
the Minkowski gauge \eqref{mk} and the minimal time function \eqref{mt} of $x=(x_1,x_2)\in\R^2$ are computed by, respectively,
\begin{equation*}
\rho_F(x)=\begin{cases}|x_1|&\mbox{ if }\;x\in\R\times\{0\},\\
\infty&\text{otherwise;}
\end{cases}
\end{equation*}
\begin{equation*}
T^F_{\O}(x)=\begin{cases}
0&\mbox{ if }\;x\in \O,\\
|x_1|-1&\mbox{ if }\;|x_2|\le 1\;\mbox{ and }\;|x_1|>1,\\
\infty&\text{otherwise.}
\end{cases}
\end{equation*}
Taking first the in-set $\ox=(1,0)\in\O$, we can easily check that $\partial
T^F_{\O}(\ox)=[0,1]\times\{0\}$ and that
$\sigma_F(v)=|v_1|$ for $v=(v_1,v_2)\in\R^2$. It is also clear that
\begin{equation*}
N(\ox;\O)\cap C^*=N(\ox;\O)\cap\big\{v\in\R^2\big|\;\sigma(-v)\le 1\big\}=[0,1]\times\{0\},
\end{equation*}
and thus we verify equality \eqref{c1} in Theorem~\ref{convex case}(i). Furthermore, it is easy
to verify that $\partial^\infty T^F_\O(\ox)=\{0\}$ in accordance with Theorem~\ref{convexsingular1}(i)
in the in-set case; this confirms that $T^F_\O$ is locally Lipschitzian around $\ox=(1,0)$

Considering another in-set point $\oy=(0,1)\in\O$, we have
\begin{equation*}
\partial T^F_{\O}(\oy)=N(\oy;\O)\cap C^*=\{0\}\times\R^+,
\end{equation*}
which verifies the conclusion of Theorem~\ref{convex case}(i). It follows similarly that $\partial^\infty
T^F_\O(\oy)=\{0\}\times[0,\infty)$, which is in accordance with Theorem~\ref{convexsingular1}(i) and
with the non-Lipschitzian behavior of the minimal time function around $\oy=(0,1)$.

Considering finally the out-of-set point $\oz=(2,1/2)\notin\O$, with the projection singleton $\Pi^F_\O(\oz)=
\{\bar w\}$ computed by $\bar w=(1,1/2)$. Then we arrive at the equalities
\begin{equation*}
\partial T^F_\O(\bar z)=\{(1,0)\}=-\partial\rho_F(\bar w-\bar z)\cap N(\bar w;\O)\;\mbox{ and }\;
\partial^\infty T^F_\O(\oz)=\{0\},
\end{equation*}
which verify the conclusions of Theorem~\ref{convexsingular1}(ii) and Theorem~\ref{convex case 2} and confirm,
in particular, the local Lipschitzian property of $T^F_\O$ around $\oz=(2,1/2)$.}
\end{Example}

\end{document}